\documentclass[10 pt, conference]{ieeeconf}

\IEEEoverridecommandlockouts                              % This command is only needed if
\usepackage{color}
\usepackage{graphics} % for pdf, bitmapped graphics files
\usepackage{needspace}

\input{mysymbol.sty}
\usepackage{theorem}
\usepackage{cite}
\usepackage{amsmath}
\usepackage{amssymb}
\usepackage{mathtools}
\usepackage{multirow}
\usepackage{graphicx}
\usepackage{url}
\usepackage{graphics, subfigure, times, amsfonts}
\usepackage{tikz, epic,eepic}
\usetikzlibrary{shapes,arrows}
\usepackage{pgfplots}
\usepackage{color}
\usepackage{hyperref}
\usepackage{epstopdf}
\usepackage{epsfig, amsbsy}
\usepackage{latexsym}
\usepackage{amscd, verbatim}
\usepackage{multirow}
\usepackage{booktabs}
\usepackage{bbm}
\usepackage{algorithm}
\usepackage{algorithmic}

\newtheorem{theorem}{Theorem}

\newtheorem{lemma}{Lemma}
\newtheorem{definition}{Definition}

{\itshape}{\rmfamily}
\newtheorem{assumption}{Assumption}
\newtheorem{remark}{Remark}
\newtheorem{condition}{Condition}

\def\forall{\text{for all\ }}

\title{\LARGE 
Decentralized Inertial Best-Response with Voluntary and Limited Communication in Random Communication Networks}
\author{Sarper Ayd\i n and Ceyhun Eksin % 
\thanks{S. Aydin and C. Eksin are with the Industrial and Systems Engineering Department, Texas A\&M University, College Station, TX 77843. E-mail:{\tt\small  \; sarper.aydin@tamu.edu; eksinc@tamu.edu}%
}}

\begin{document}
\normalsize
\maketitle

%%%%%%%%%%%%%%%%%%%%%%%%%%%%%%%%%%%%%%%%%%%%%%%%%%%%%%%%%%%%%%%%%%%%%%%%%%%
%%%   A   B   S   T   R   A   C   T   %%%%%%%%%%%%%%%%%%%%%%%%%%%%%%%%%%%%%
%%%%%%%%%%%%%%%%%%%%%%%%%%%%%%%%%%%%%%%%%%%%%%%%%%%%%%%%%%%%%%%%%%%%%%%%%%%
%
\begin{abstract}
Multiple autonomous agents interact over a random communication network to maximize their individual utility functions which depend on the actions of other agents. We consider decentralized best-response with inertia type algorithms in which agents form beliefs about the future actions of other players based on local information, and take an action that maximizes their expected utility computed with respect to these beliefs or continue to take their previous action. We show convergence of  these types of algorithms to a Nash equilibrium in weakly acyclic games under the condition that the belief update and information exchange protocols successfully learn the actions of other players with positive probability in finite time given a static environment, i.e., when other agents' actions do not change. We design a decentralized fictitious play algorithm with voluntary and limited communication (DFP-VL) protocols that satisfy this condition. In the  voluntary communication protocol, each agent decides whom to exchange information with by assessing the novelty of its information and the potential effect of its information on others' assessments of their utility functions. The limited communication protocol entails agents sending only their most frequent action to agents that they decide to communicate with. Numerical experiments on a target assignment game demonstrate that the voluntary and limited communication protocol can more than halve the number of communication attempts while retaining the same convergence rate as DFP in which agents constantly attempt to communicate.
\end{abstract}

%%%%%%%%%%%%%%%%%%%%%%%%%%%%%%%%%%%%%%%%%%%%%%%%%%%%%%%%%%%%%%%%%%%%%%%%%%%
%%%   E N D     A   B   S   T   R   A   C   T   %%%%%%%%%%%%%%%%%%%%%%%%%%%%%%%%%%%%%
%%%%%%%%%%%%%%%%%%%%%%%%%%%%%%%%%%%%%%%%%%%%%%%%%%%%%%%%%%%%%%%%%%%%%%%%%%%
%

%%%%%%%%%%%%%%%%%%%%%%%%%%%%%%%%%%%%%%%%%%%%%%%%%%%%%%%%%%%%%%%%%%%%%%%%%%%
%%%   S E C T I O N %%%%%%%%%%%%%%%%%%%%%%%%%%%%%%%%%%%%%
%%%%%%%%%%%%%%%%%%%%%%%%%%%%%%%%%%%%%%%%%%%%%%%%%%%%%%%%%%%%%%%%%%%%%%%%%%%
%
\section{Introduction}

Multi-agent systems comprise of interlinked decision-makers (agents) aiming to maximize  objectives that depend on the actions of other agents in the system. In epidemics, the preemptive measures taken by individuals affect the risks associated with socialization \cite{eksin2017disease,bauch2004vaccination}. In a smart grid, multiple devices determine generation and consumption levels to reach a  balance while minimizing costs \cite{kar2014distributed,zhang2012robust}. In autonomous teams of mobile robots, each robot decides its direction of movement and position to maximize a team objective that depends on the movements and positions of other robots \cite{aydin2020communication,kantaros2016distributed,kantaros2019temporal}. In all of these settings, agents have to reason about the motives of other  agents based on local information. Game theoretic equilibrium concepts, i.e., Nash equilibrium (NE), provide a benchmark for rational reasoning where agents assume other agents are also trying to maximize their individual or team objectives. However, computation of NE is not feasible given limited computation capabilities and local information. Here, we develop decentralized game-theoretic learning algorithms for settings where agents do not know the incentives of other agents, and need to communicate over a random network that is subject to failures in order to reason about other agents' actions.

Success of a communication attempt is often subject to random failures in social and technological settings. Moreover, in social settings communication is often voluntary, i.e., agents attempt to communicate upon the need for information exchange. In technological settings, communication is costly to the agents. Because of this, decentralized learning algorithms in which agents constantly attempt to communicate are neither realistic representations of information exchange in social settings, nor practical in technological settings. Here, we propose decentralized learning algorithms in which agents consider the effect of their information on a potentially receiving agent's beliefs before attempting to communicate.

% Communication often occurs in random communication networks that is subject to failures... Not only that but communication is also costly---> there is a need to develop communication efficient learning algorithms
In the decentralized game-theoretic learning algorithms considered in this paper, agents use best-response with inertia to determine their next actions at each step. In best-response with inertia, each agent forms beliefs about the actions of other agents, and takes an action that either maximizes its expected utility computed with respect to its beliefs (best-responds) or continues to take its former action (shows inertia). Whether an agent best-responds or shows inertia in a given step is random. Agents form beliefs about other agents' behavior via information exchanges over a random communication network. The randomness of communication means that agents cannot receive information from every other agent at each step. Given this setting and learning updates, we show convergence of the best-response with inertia behavior to a NE of any weakly acyclic game in finite time almost surely, as long as the information exchange and belief update protocols ensure that agents are able to learn another agent's action if that agent repeats the same action long enough (Theorem \ref{thm_main}). 

We call this sufficient condition for convergence (Condition \ref{cond_prob}) as {\it prediction under static actions}. Based on this condition, we design voluntary communication protocols in which agents attempt to send information to an agent if they see the need to communicate (Section \ref{sec_com}). Agents determine the need to communicate based upon the novelty of their information to the potential receiving agent. That is, each agent assesses the novelty of their recent information to other agents. For such an assessment, agents form second order beliefs, i.e., reason about the beliefs that other agents have about their behavior. In this voluntary communication protocol, agents assume other agents act according to a stationary distribution determined by the past empirical frequencies of their actions similar to standard fictitious play (FP) \cite{brown1951iterative,young2004strategic,marden2009joint}. Unlike FP, agents cannot keep track of the empirical frequencies of all the agents when the communication is random and voluntary.
We show that the voluntary communication protocol satisfies the prediction under static actions condition when agents attempt to send only the frequencies of their most frequent actions (Theorem \ref{lem_prob_com2}). Via numerical experiments in a target assignment problem, we demonstrate that the proposed DFP algorithm with voluntary communication and limited information exchange (DFP-VL) can perform similar to DFP with constant communication attempts in convergence rate while more than halving the communication attempts per link (Section \ref{sec_numeric}).

% \red{We note that the agents do not have to repeat the same action, the condition 

% We assume agents form beliefs about other agents' actions based on local information exchanges occurring over random communication links. }

\subsection{Related Literature}

FP converges to rational behavior in various games including potential \cite{monderer1996fictitious}, weakly acyclic \cite{young2004strategic,marden2009cooperative}, zero-sum \cite{robinson1951iterative}, near-potential \cite{candogan2013dynamics}, and stochastic games \cite{sayin2020fictitious}. In FP, each agent takes an action that maximizes its expected utility (best responds) assuming other agents select their actions randomly from a stationary distribution. Agents assume this stationary distribution is given by the past empirical frequency of past actions. FP is not a decentralized algorithm, since agents need to observe past actions of everyone to be able to keep track of empirical frequencies, and compute the expectations of their utilities. Recent works \cite{Swenson_et_al_2014,eksin2017distributed,swenson2018distributed,arefizadeh2019distributed} consider a decentralized form of the fictitious play, in which agents form estimates on empirical frequencies of other agents' actions by averaging the estimates received from their neighbors in a communication network. These algorithms are shown to converge to a NE in weakly acyclic games, i.e., games that admit finite best-response improvement paths. However, they rely on communication with neighbors after every decision-making step. This assumption ignores the randomness of communication attempts, e.g., in wireless communication settings, and the energy costs of communication. Preliminary versions of this paper either consider a specific setting for the voluntary communication protocol design, namely the target assignment game in \cite{aydin2020communication}, or focus on the convergence of a specific communication protocol for DFP in \cite{aydin2020decentralized}.
Theorem \ref{thm_main} generalizes the convergence results in DFP and the preliminary results in \cite{aydin2020communication,aydin2020decentralized} to show that a generic inertial best-response type behavior will converge to a rational action profile as long as there exists a belief update and information exchange protocol in which agents are able to learn the actions of other agents when these agents repeat the same action long enough. We then leverage this result to design an intuitive and novel class of communication efficient belief update and information exchange protocols.

% In this paper, we design a decentralized communication protocol for the DFP that allows agents to determine whom to communicate with and when to cease communication (Section \ref{sec:model}). The communication protocol is based on the idea that agents do not need to send their empirical frequencies to other agents unless they carry new information that can potentially affect the decisions of other agents. In the context of DFP, we realize this idea by linking novelty of information and potential effect on others' evaluations to two metrics computed respectively as the change in the empirical frequency caused by the current action, and the error others make in estimating the agent's empirical frequency. After taking a new action, each agent attempts to transmit its empirical frequency to another agent if the values of these two metrics exceed certain constant thresholds. Our main result (Theorem \ref{thm_main}) shows that the DFP-V algorithm with the threshold-based communication protocol converges to a pure NE action profile in finite time given small enough thresholds, if agents' beliefs about the state of the environment weakly converges to a common belief. Numerical experiments demonstrate that the communication protocol can lower the communication attempts by half while showing similar convergence properties as the standard DFP algorithm (Section \ref{sec_numeric}). 

In the voluntary information exchange protocols, the assessment of the novelty of information to a potential receiving agent is based on two metrics: {\it i)} novelty of local information and {\it ii)} its potential effect on the belief of the receiving agent. Such metrics that are based on second order beliefs (estimating the estimates of the receiving agents) can provide similar benefits to communication efficiency in other decentralized game-theoretic learning algorithms based on, e.g., gradient descent \cite{alpcan2005distributed,koshal2016distributed,Shamma_Arslan_2005,de2019distributed}, best-response \cite{scutari2013joint}, ADMM \cite{salehisadaghiani2019distributed}, and other adaptive strategies \cite{ye2021adaptive}. Indeed, communication-censoring based protocols that rely on some form of novelty of information metrics recently proved viable in reducing communication attempts in distributed stochastic gradient descent \cite{chen2018ordered,chen2018lag} and ADMM \cite{li2019communication} in the context of optimization. In the class of information exchange protocols considered here, while the novelty of information metric is sender specific, the metric on potential effect of information on other's assessment is receiving agent specific. Thus, agents manage their local information by deciding whom to communicate with. This is a novel communication protocol that relies on agents keeping track of second order beliefs, i.e., forming beliefs on beliefs, in order to estimate the novelty of their information to the candidate receiving agent. 
\section{Multi-Agent Systems in Time-Varying Random Networks}\label{sec:model}

\subsection{Notation}

We use $||.||$ to denote Euclidean norm. The notation $\Delta(.)$ defines the space of probability distributions over given set. $\bbone_{(.)}$ is indicator function.  Its value is $1$ if the given condition is satisfied, otherwise $0$.
\subsection{Game-Theoretical Definition}
We consider a strategic game $\Gamma$ among a set of $N$ agents denoted with $\ccalN=\{1,2,\cdots, N\}$. Each agent $i$ chooses an action $a_i$ from a common action set $\ccalA$ with finitely many actions, i.e., $|\ccalA|=K$. We represent each action with an unit vector $\bbe_k \in \reals^K$ so that  $\ccalA:=\{\bbe_1,\bbe_2,\cdots,\bbe_K\}$. Each agent has an utility function $u_i: \ccalA^N \rightarrow \reals$ that depends on the joint action profile $(a_i,a_{-i})\in\ccalA^N$ where $-i$ denotes the set of all agents, and $a_{-i}$ is the action profile of agents in the set $-i$. The strategic game $\Gamma$ is defined by the tuple $(\ccalN,\mathcal{A}^N,\{u_i\}_{i\in\ccalN})$. 

A mixed action (strategy) $\sigma_i$ is a probability distribution over the action space. We define the space of probability distributions over the action space as $\Delta (\ccalA)$. A strategy profile $\sigma=(\sigma_i, \sigma_{-i})$ is a joint mixed action profile belonging to the set of independent probability distributions over the space of action profiles, i.e., $\Delta^N(\ccalA) = \prod_{i \in \ccalN} \Delta(\ccalA)$. We denote the expected utility of agent $i$ given a strategy profile $\sigma$ as $u_i(\sigma) := \sum_{a\in\ccalA^N} u_i(a)\sigma(a)$ where $\sigma(a)$ is the probability of action profile $a \in \ccalA^N$.
Next, we describe the standard FP and then introduce a generalization of FP for random communication networks. 
% Next we introduce the standard FP algorithm. 

%%%%%%%%%%%%%%%%%%%%%%%%%%%%%%%%%%%%%%%%%%
\subsection{Fictitious Play with Inertia} 
%\grey{We show the selection of agent $i$ at time $t\in \naturals^+$, by $a_i(t)\in \ccalA$.}
%In (centralized) FP algorithm, it is assumed that agents repeatedly select actions  according to a stationary distribution that is determined by the histogram of their past actions. The histogram, i.e., the empirical frequency $f_i(t)$ of agent $i$ is as follows,

%\red{Fictitious Play (FP) is a game-theoretical learning algorithm where agents are assumed to take actions from a stationary probability distribution.}
FP is a distributed game-theoretic learning algorithm in which agents repeatedly take actions in discrete time steps $t=1,2,\dots$ that maximize their expected utilities computed with respect to some estimate of the other agents' strategies. In estimating the strategies of others, each agent assumes other agents are taking actions drawn from a stationary probability distribution determined by the {\it empirical frequency} of the past actions of agents. The empirical frequency of agent $i$ is computed as follows,
\begin{equation} \label{eq_empirical_frequency}
    f_i(t) = (1-\rho)f_i(t-1)+\rho a_i(t),
\end{equation}
where $a_i(t)\in \ccalA$ is the action of agent $i$ at time $t\in \naturals^+$ and $\rho \in(0,1)$ is a fading memory constant determining the update rate of the empirical frequency. 

%\red{Therefore, agents keep statistics of past actions.} \green{CE: Avoid unnecessary conjunctions when describing models.} \red{These statistics are named as \textit{empirical frequency} and computed as below,}
%
%\red{Then, using the fact that $f_i(t) \in \Delta \ccalA$, we can explicitly define utility functions as the expected value with respect to empirical frequencies, }
Given the empirical frequencies of other agents $f_{-i}(t)=\{f_j(t)\}_{j\in \ccalN\setminus i}$, agent $i$'s expected utility from taking action $a_i$ is given as, 
\begin{equation} \label{eq_exp_utility}
u_i(a_i, f_{-i}(t)) = \hspace{-4pt}\sum_{a_{-i} \in \ccalA^{N-1}}  u_i(a_i, a_{-i})  f_{-i}(t)(a_{-i}).
\end{equation}

In FP with inertia, each agent best-responds with inertia, i.e., either takes an action that maximizes its expected utility, or follows its previous action with a small probability $\epsilon\in(0,1)$. 
%
% \begin{align}\label{eq_response}
% a_{i}(t)=
% \begin{cases}
% \argmax_{a_i\in\ccalA} u_i(a_i,f_{-i}(t-1)) &\;\; \text{w.pr. } 1-\epsilon ,\\
% a_{i}(t-1) &\;\; \text{w.pr. } \epsilon.
% \end{cases}
% \end{align} 
Agent $i$ needs to observe the past actions of {\it all} agents in order to compute the empirical frequencies as per \eqref{eq_empirical_frequency} so that it can compute the best response action.

\subsection{Decentralized fictitious play (DFP) in random networks}

When communication between agents is subject to failures, agents do not have immediate and permanent access to others' actions. One way to address this problem is by agents keeping local estimates of empirical frequencies of past actions, in which each agent forms estimates about the empirical frequencies of other agents based on information received from neighboring agents in the communication network. 

The estimate of agent $i$ on agent $j$'s empirical frequency in \eqref{eq_empirical_frequency} is denoted with $f^i_j(t)\in\Delta(\ccalA)$. Replacing the empirical frequencies $f_{-i}(t)$ with the estimates $f^i_{-i}(t):=\{f^i_j(t)\}_{j\in\ccalN\setminus\{i\}}$ in \eqref{eq_exp_utility}, we get the expected utility of agent $i$ from taking action $a_i \in\ccalA$ denoted as $u_i(a_i, f^i_{-i}(t))$. As in standard FP with inertia, agents best-respond with inertia, i.e., maximize their expected utility or continue taking the previous action,
\begin{align}\label{eq_response_dc}
a_{i}(t)=
\begin{cases}
\argmax_{a_i\in\ccalA} u_i(a_i,f^i_{-i}(t-1)) &\;\; \text{w.pr. } 1-\epsilon ,\\
a_{i}(t-1) &\;\; \text{w.pr. } \epsilon.
\end{cases}
\end{align} 

In DFP, agents update their local estimates based on information they receive from agents that send information over the random communication network. 

Specifically, we assume point-to-point communication between each pair of agents is possible but communication is subject to random failures. The probability of the existence of a communication link between agent $i \in \ccalN$ and agent $j \in \ccalN \setminus \{i\}$ at time $t\in \naturals^+$ is distributed with a Bernoulli random variable, 
\begin{align} \label{eq_prob}
 c_{ij}(t)  \, \sim \text{Bernoulli} (p_{ij}(t)),
\end{align}
where the {probability of success is $0 \le  p_{ij}(t)<1$}. 

We denote the random communication network at time $t$ with $G(t)=(\ccalN, \ccalE(t))$ where $\ccalE(t)$ is the set of edges realized according to \eqref{eq_prob}. The random communication network $G(t)$ belongs to the space of all possible networks $\ccalG$.

We denote the history of the actions and networks up to time $t$ as {$H(t):=(\ccalA^N \times \ccalG)^t$}. We define a measurable space {$(H(\infty),\ccalB)$} as the sequence of actions and networks $H(\infty)$ and the Borel sigma-algebra ($\ccalB$). We let {$\{\ccalH(t)\}_{t \ge 0}$} be a sub-sigma algebra of $\ccalB$. The information available to agent $i$ at time $t$ is denoted with ${\ccalH_{i}(t)}$.

The information exchange protocol of agent $i$, denoted with {$\Omega_i: H(t) \to \ccalN\times H(t)$}, determines the set of agents agent $i$ is willing to communicate with (${\ccalN_{i}^{out}(t)\subseteq\ccalN}$) and the information agent $i$ shares with them that is measurable with respect to the information available ${\ccalH_{i}(t)}$. Upon receiving information from its neighbors ${\ccalN_{i}^{in}(t):=\{j\in\ccalN\setminus\{i\}: i\in \ccalN_j^{out}(t)\}}$, agent $i$ updates its estimates about the empirical frequencies of other agents $\{f^{i}_{j}(t)\}_{j\in\ccalN}$ according to a function {$\Phi_{i,j}:H(t)\to \Delta(\ccalA)$} that is measurable with respect to the information available at time $t+1$ ({$\ccalH_{i}(t+1)$}). We let $H_{i}(t)$ be the realization of the information available to agent $i$, $H_{i}(t):=\{\{a_{i}(s)\}_{s=1}^{t-1},\prod_{s=1}^{t-1}\prod_{j\in\ccalN_i^{in}(s)} \Omega_j(H_{j}(s))\}$ as a result of the exchange protocol $\{\Omega_j\}_{j\in\ccalN}$. The exchange protocol determines the information available to each agent in the next time step $H_{i}(t+1) = \{H_{i}(t), a_i(t), \prod_{j\in\ccalN_i^{in}(t)} \Omega_j(H_{j}(t))\}$. For the convergence analysis, we will be agnostic to the specifics of the estimate updates $(\Phi_{i})$  and the information exchange process $(\Omega_i)$, as long as they ensure that agents are able to learn others' actions under a static action profile. We state the condition formally next. 

\begin{condition}[Prediction under static actions]\label{cond_prob}
There exists a positive probability $\hat{\epsilon}>0$ and a finite time $\hat T$ such that if an agent $j\in\ccalN$ repeats the same action for at least $T>\hat T$ times starting from time $t>0$, i.e., $a_{j}(s)=\bbe_k$ for  $s=t,t+1,\cdots,{t+T-1}$ and $\bbe_k\in\ccalA$, agent $i\in\ccalN$ learns agent $j$'s action with positive  probability $\hat{\epsilon}>0$, i.e., $\mathbb{P} (||a_j(t+T)-f^i_j(t+T)|| \le {\xi}| \ccalH(t)) \ge \hat{\epsilon}$ for any $\xi>0$. %\green{CE: USE another greek letter here. $\eta$ is used for parameters of the DFP-LV...}
\end{condition}

Any estimate update and information exchange process that satisfies Condition \ref{cond_prob} makes sure that agent $i$'s estimate of agent $j$'s action  $f^i_j(t)$ gets close to agent $j$'s action whenever agent $j$ repeats its action long enough. 

We summarize key steps of the generic DFP  next. 
\begin{algorithm}[H] 
   \caption{Generic DFP for Agent $i$}
\label{suboptimal_alg_DFP}
\begin{algorithmic}[1]
   \STATE {\bfseries Input:} Inertia probability $\epsilon$ and fading constant $\rho$.
   \STATE {\bfseries Given:} $f_{-i}^i(0)$ and $a(0)$ for all $i \in \ccalN$.
\FOR{$t=1,2,\cdots $}
  \STATE {\it Best-respond}: Use $f^i(t-1):=\{f^i_j(t-1)\}_{j\in\ccalN\setminus i}$ in \eqref{eq_response_dc} \STATE {\it Share information}: Use $\Omega_i$ to determine ${\ccalN_{i}^{out}(t)}$ and information to be exchanged
\STATE {\it Observe}: Receive information from ${\ccalN_{i}^{in}(t)}\cap \{j: c_{ji}(t)=1\}$ 
\STATE {\it Update estimates}: $f_j^i(t+1)= {\Phi_{i,j}(H_{i}(t))}$ for $j\in\ccalN\setminus i$.
  \ENDFOR 
   \end{algorithmic}
\end{algorithm}

%\green{\large CE:  this assumption should go to the next section...} 
%\blue{We assume the success probability is bounded by some positive constant $0<\rho_c<p_{ij}(t)$ for all $t\in \naturals^+$ and $i,j\in\ccalN$.}

%%%%%%%%%%%%%%%%%%%%%%%%%%%%%%%%%%%%%%%%%%%%%%%%%%%%%%%%%%%%%%%%%%%%%%%%%%%
%%%   S E C T I O N %%%%%%%%%%%%%%%%%%%%%%%%%%%%%%%%%%%%%
%%%%%%%%%%%%%%%%%%%%%%%%%%%%%%%%%%%%%%%%%%%%%%%%%%%%%%%%%%%%%%%%%%%%%%%%%%%
%

%%%%%%%%%%%%%%%%%%%%%%%%%%%%%%%%%%%%%%%%%%%%%%%%%%%%%%%%%%%%%%%%%%%%%%%%%%%
%%%   S E C T I O N %%%%%%%%%%%%%%%%%%%%%%%%%%%%%%%%%%%%%
%%%%%%%%%%%%%%%%%%%%%%%%%%%%%%%%%%%%%%%%%%%%%%%%%%%%%%%%%%%%%%%%%%%%%%%%%%%

\section{DFP Convergence for Weakly Acyclic Games} \label{sec_con}
%\green{CE: This section needs to be updated according to the updates in the previous section}

 %\grey{DFP-V (Algorithm \ref{suboptimal_alg}) is structured in the  decentralized form that agents has autonomous control over both its action selections, and when and whom to exchange information with. In this section, we analyze the convergence of action profile, that no agent wish to change its strategy. Before delving into the details, we introduce the concepts that the convergence proofs will employ. }
%\green{CE: too many sentences with ``that". Most of them have no place in this introductory paragraph. Form simple sentences.}

%Here we show convergence of the action profiles under DFP-V to a pure NE of the game in finite time for particular class of games, called weakly acyclic games \cite{young1993evolution,milchtaich1996congestion}. Below, we define weakly acyclic games.

%\red{The convergence theorem are going to be based on game-theoretical concepts. we begin with required definitions and assumptions related with game theory to show the convergence. } \green{CE: Remove. Carries no information. You do not need to prep the reader.}

We consider convergence of the DFP in the class of weakly acyclic games which have (finite) sequence of best-response updates that end up at a pure Nash equilibrium, named as finite improvement paths \cite{young1993evolution,milchtaich1996congestion}.

A Nash equilibrium strategy is an (mixed) action profile in which no individual agent can benefit by unilaterally switching to another action. A formal definition follows.

\begin{definition} [Nash Equilibrium] \label{def_Nash}
A strategy profile  $\sigma^*=(\sigma_i^*,\sigma_{-i}^*) \in \Delta^N(\mathcal{A})$ is a Nash equilibrium of the game $\Gamma$ if and only if for all $i\in\ccalN$
\begin{equation} \label{eq_Nash}
    u_i(\sigma^*_i,\sigma^*_{-i})\ge u_i(\sigma_i,\sigma_{-i}), \quad \forall \sigma_i \in \Delta(\mathcal{A}).
\end{equation}
A pure NE strategy profile $\sigma^*$ is a NE that selects an action profile $a=(a_i,a_{-i}) \in \ccalA^N$ with probability 1. 
\end{definition}

A best-response path is a sequence of action profiles obtained by a single agent best-responding to the current action profile at each step of the sequence. Next, we provide a formal definition of weakly acyclic games. 

\begin{definition} [Weakly Acyclic Games] \label{def_WA}
A game $\Gamma$ is weakly acyclic if from any joint action profile $a=(a_i,a_{-i})\in\ccalA^N$, there exists a best-response path ending at a pure NE $a^*=(a^*_i,a^*_{-i})$.
\end{definition}

The existence of a finite best-response path ensures that no agent can improve its utility after some finite number of iterations. 
Weakly acyclic games are {a} broad class of games that include potential games and its several variants such as best-response potential and pseudo-potential games. %The analysis will utilize the existence of finite improvement path to show almost sure converge to a pure NE.

We consider weakly acyclic games in which optimal action is unique against others' actions if other agents take actions according to a pure NE action profile. Specifically, we make the following assumption.

\begin{assumption}\label{as_single}
For any pure NE action profile $a^*\in \ccalA^N$ of the game $\Gamma$, it holds that,
\begin{equation}
  \{a^*_i\} = \argmax_{a_i\in\ccalA} u_i(a_i,a^*_{-i}).
\end{equation}
\end{assumption}
This assumption makes sure that agents are not indifferent between multiple actions at a pure Nash equilibrium.

\subsection{Convergence to a Pure Nash Equilibrium}

We show almost sure convergence of joint action profile $a(t)$ to a pure NE $a^*$ (Theorem \ref{thm_main}). 
The convergence result relies on the fact that action profile stays forever at a pure NE once it reaches the NE (Lemma \ref{lem_ab}), and there is a positive probability to reach a pure NE from any action profile(Lemma \ref{lem_pos_ab}). Before showing these lemmas, we show that the best response action of an agent computed with respect to the estimated empirical frequencies  $\{f^i_j(t)\}_{j\in \ccalN}$ belongs to the best response action set computed with respect to the actual actions of others $a_{-i}(t)$, whenever the estimates are close enough to $a_{-i}(t)$--see Appendix  \ref{ap_pr_2} for the proof.

%\red{start with proving that agents take the optimal action set against local estimates of others' empirical frequencies $f^i_j(s)$, as they learn others' actions $a_{-i}$.  }\green{CE: Remove}

\begin{lemma}\label{lem_est}
There exists a small enough ${{\xi>0}}$ such that if $||a_j(t)-f^i_j(t)|| \le {\xi}$ for agents $j \in \ccalN \setminus \{i\}$ at time step $t$, then $\argmax_{a_i \in \ccalA} u_i(a_i,f^i_{-i}(t)) \subseteq \argmax_{a_i \in \ccalA} u_i(a_i,a_{-i})$ for all $i\in\ccalN$.
\end{lemma}
%\green{CE:Why do you need to write this for $t+T$? For instance, for anytime $t$ there exists a small enough $\eta$ such that if $||a_j(t)-f^i_j(t)|| \le \eta$, then ...}
% \begin{myproof}
% See Appendix \ref{ap_pr_2}.
% \end{myproof}

Next, we prove that pure NE have absorption property. When agents play a pure NE and are aware of others' actions, agent are going to stay in this pure NE indefinitely.

\begin{comment}

\begin{lemma}\label{lem_prob_repeat}\textbf{(positive probability of repetition and communication)}
Suppose Assumption \ref{as_filt} holds and condition in \eqref{eq_condition} is not true  $\forall{(i,j)} \in \ccalN \times \ccalN \setminus \{j\}$. Let $E_1$ be the event is defined follows ,
\begin{align} \label{eq_event}
    E_1(t)=\{&a(s)=a,c_{ij}(t+T)=1, \nonumber\\
    &\forall{s} \in \{t,t+1,\cdots,t+T\}\},
\end{align}
where $a(s)$ is a joint action profile at time $s$ and $c_{ij}(t)$ is the realization of Bernoulli random variable determining communication link between $i$ and $j$. Then, the probability of the event $E_1(t)$ conditioned on $\ccalF(t)$ is bounded below by a positive constant $\epsilon_1(T)$,
\begin{equation}
    \mathbb{P}(E_1(t)|\ccalF_t,\beta_{ij}(t+T)=1) \ge \epsilon_1(T).
\end{equation}
\end{lemma}
\begin{myproof}
The proof follows from the fact that due to inertia, the probability of repetition of actions in finite $T$ time is always positive. Further, there is a positive probability of communication (at least $p_{ij}(t)>\epsilon_{com}$) before condition \eqref{eq_condition} is satisfied. Thus the probability of action repetitions and communication is at least $\epsilon_1=\epsilon^{NT}\epsilon_{com}^{N(N-1)}$. 
\end{myproof}
\end{comment}

\begin{lemma}\label{lem_ab} \textbf{(absorption property)}
Suppose Assumption \ref{as_single} holds. Assume $||a_j(t+T)-f^i_j(t+T)|| \le {\xi}$ where $\xi>0$ satisfies Lemma \ref{lem_est}  for all pairs of agents $(i,j) \in \ccalN \times \ccalN \setminus \{i\}$ at time step $t+T$. Further, let $a^*\in \ccalA^N$ be a pure NE action profile and $a(t+T)=a^*$. Then, $a(s)=a^*=(a_1^*,a_2^*,\cdots,a_N^*)$ holds, $\forall{s} \ge t+T$.
\end{lemma}
 
\begin{myproof}
% By Lemma \ref{lem_est}, it holds, $\argmax_{a_i \in \ccalA} u_i(a_i(t+T),f^i_{-i}(t+T)) \subseteq \argmax_{a_i \in \ccalA} u_i(a_i,a_{-i}(t+T))$. Since, $a(t+T)=a^*$ and b
By Assumption \ref{as_single}, the set of optimal actions given others' actions $a_{-i}(t+T)=a^*_{-i}$ is a singleton given by $\argmax_{a_i \in \ccalA} u_i(a_i,f^i_{-i}(t+T)) = \argmax_{a_i \in \ccalA} u_i(a_i,a^*_{-i})=\{a_i^*\}$. Otherwise, by inertia agent $i$ takes the the same action $a_i^*$. Thus, the joint action profile remains at the pure NE, i.e., $a(s)=a^*, \, \forall{s} \ge t+T$.
\begin{comment}
Then, by definition of a pure NE $\{a_i^*\}\in\argmax_{a_i \in \ccalA} u_i(a_i,f^i_{-i}(s))$ for $s\ge t+T_1+T_2$. Moreover, by Lemma \ref{lem_single}, the set $\argmax_{a_i \in \ccalA} u_i(a_i,f^i_{-i}(s))$ is a singleton for each $i \in\ccalN$. Thus the action profile $a^*$ is repeated at time $t+T_1+T_2$. Inductively, $a(s) = a^*$ for all $s \geq t$.
\end{comment}
\end{myproof}

The next lemma states that there is a positive probability that agents can reach a NE action profile with any communication scheme that satisfies Condition \ref{cond_prob}.

\begin{lemma}[\textbf{positive probability of absorption}]\label{lem_pos_ab}
Suppose Assumption \ref{as_single} and Condition \ref{cond_prob} hold. Let $a(t)$ be the joint action profile at time $t$ and $f^i(t):=\{f^i_j(t)\}_{j\in\ccalN}$ be agent $i$'s estimate on all agents at time $t$. At time $t$, we define the following event $\forall{(i,j)} \in \ccalN \times \ccalN \setminus \{i\}$, 
\begin{align*}
     E(t)=&\{a(s)=a^*, \, ||a_j(\bar{s}+T)-f^i_j(\bar{s}+T)||  \le {\xi} \\
     %&\text{for some}\, m \in \{t, 1, \cdots,(K^N+1)T\}\}
     &\forall{s} \in \{\bar{s},\bar{s}+1,\cdots,\bar{s}+T-1\} \\
     %&\forall{s_2} \in \{\bar{s}+{T_1}, \bar{s}+{T_1}+1,\cdots,\bar{s}+{T_1}+{T_2}-1\} \\
     &\text{for some}\, \bar{s} \in \{t, t+1, \cdots,
     t+K^NT\}\}
\end{align*}
where $a^*$ is a pure NE. There exists ${\xi}>0$  small enough such that the transition probability $\mathbb{P}({E(t)}|{\ccalH(t)}) \ge \bar{\epsilon}(T)$, is bounded below by $\bar{\epsilon}(T)>0$ and always positive $\forall{t} \in \naturals^+$.
\end{lemma}
\begin{myproof}
To show the result, we are going to use the fact that in weakly acyclic games, there exists a finite path from any action profile to a pure NE. Since, the action  set of each agent is finite and its size is equal to $K$, there exist $K^N$ different joint action profiles in total. Hence, it is the upper bound on the length of finite path to a pure NE. Thus, if $a(t)=a^*$, the pure NE is reached and, the proof is trivially completed.

If $a(t)\not =a^*$, we are going to exploit the fact that the finite path to a pure NE consists of finite improvement paths. In each improvement path, only one agent improve its utility by changing its action. Therefore, all agents firstly stay in their actions  
so that it holds  $||a_j(t+T)-f^i_j(t+T)|| \le {\xi}$ with probability at least $\hat{\epsilon}^N$ by Condition \ref{cond_prob}. Then, by inertia, at time step $t+T$, there exists a positive probability $\epsilon^{(N-1)}$, $N-1$ agents continue to stay in the same action for one more time step, and only one agent takes the optimal action against others with probability $(1-\epsilon)$. Since it holds $\argmax_{a_i \in \ccalA} u_i(a_i,f^i_{-i}(t+T)) \subseteq \argmax_{a_i \in \ccalA} u_i(a_i,a_{-i})$ by Lemma \ref{lem_est}, a finite improvement path can complete with at least the probability $\epsilon_1=\hat{\epsilon}^N(1-\epsilon)\epsilon^{(N-1)}$.

After the completion of an improvement path, the event of another improvement path until $a^*$ is reached has at least the same positive probability $\epsilon_1$. As stated before, total number of improvement paths cannot exceed $K^N$ times. Once $a(\bar{s})=a^*$, the probability of learning other's actions is again $\hat{\epsilon}^N$ corresponding to the event that all agents repeat their actions at least $T$ times. Using this, the probability to reach a pure NE is bounded below as $\mathbb{P}({E(t)}|{\ccalH(t)}) \ge \bar{\epsilon}= \epsilon_1^{K^N}\hat{\epsilon}^N$. 
\end{myproof}
%

%Lemma \ref{lem_pos_ab} states that DFP-V can follow finite-best response path with positive probability. To show this,  Lemma \ref{lem_pos_ab} utilizes property of inertia and the communication protocol defined by \eqref{eq_condition}. Finally,  we  state our main convergence theorem.  

Above result leverages the fact that as long as each agent recognizes others' empirical frequencies converge to a pure action profile when they continue to take the same action, an agent can improve its utility. Now, we are ready to state the main convergence theorem.  

\begin{theorem}\label{thm_main}
Suppose Assumption \ref{as_single} and Condition \ref{cond_prob} hold. Let $\{a(t)=(a_1(t),(a_2(t),\cdots,a_N(t))\}_{t\ge1}$ be a sequence of actions by the DFP Algorithm (Algorithm \ref{suboptimal_alg_DFP}) and random time-varying communication networks $\{G(t)\}_{t\ge 1}$. The action sequence $\{a(t)\}_{t\ge1}$ converges to a pure NE $a^*$ of the game $\Gamma$, almost surely. 
\end{theorem}
%\green{CE: You need the replace this with the sequence of networks here?} 
\begin{myproof}
By Lemma \ref{lem_ab}, pure Nash equilibria are the only absorbing states among joint action profiles. By Lemma \ref{lem_pos_ab}, there exists a positive probability to reach a pure NE. Therefore, in finite time with probability $1$, a pure NE is reached and action profile stays same once reached. Thus,  the action sequence $\{a(t)\}_{t\ge1}$ converges to a pure NE $a^*$ of the game $\Gamma$, almost surely. 
\end{myproof}

The convergence theorem relies on the idea of absorbing Markov chains in which pure Nash equilibria are the only absorbing states among all joint action profiles (states). We proved almost sure convergence of actions to a pure NE by the existence of finite improvement paths and the fact that reaching a pure NE from any joint action profile has a positive probability. 
%\red{In the next section, we provide more detailed and specific communication protocols that aims to reduce communication overhead. }
%\green{CE: No telegraphing what you will do. Here you can provide another perspective on this result saying that we showed that the action profile sequences is a Markov Chain with an absorbing states as the Nash equilibria action profiles.Also add a remark discussing how this result generalizes prior similar results in our papers.}

%%%%%%%%%%%%%%%%%%%%%%%%%%%%%%%%%%%%%%%%%%%%%%%%%%%%%%%%%%%%%%%%%%%%%%%%%%%
%%%   S E C T I O N %%%%%%%%%%%%%%%%%%%%%%%%%%%%%%%%%%%%%
%%%%%%%%%%%%%%%%%%%%%%%%%%%%%%%%%%%%%%%%%%%%%%%%%%%%%%%%%%%%%%%%%%%%%%%%%%%
\section{Information Exchange and Belief Update Protocols for Random Communication Networks} \label{sec_com}

%\red{Throughout this section, we aim to introduce local communication protocols in order to reduce communication costs, while ensuring theoretical convergence.  Our discussion is going to include the comparison with DFP outlined as Algorithm \ref{suboptimal_alg_DFP} in terms of algorithmic complexities. We further are going to use additional mathematical definitions. when required, and we start with intermittent communication. } 

We introduce information exchange $\Omega_i(\cdot)$ and belief update $\Phi_i(\cdot)$ protocols that aim to reduce the number of communication attempts while at the same time guaranteeing that prediction under static actions condition (Condition \ref{cond_prob}) holds. 

\subsection{Voluntary Communication Protocols} 
%
% Agent $i$ attempts to send its empirical frequency $f_i(t)$ if it decides to communicate to agent $j$. 
We use two metrics, novelty and belief similarity, to determine whether agent $i$ attempts to communicate to agent $j$ or not. The novelty metric is the distance between the empirical frequency of agent $i$ and its current action denoted with $h_{ii}(t):= ||a_i(t)- f_i(t)||$. The belief similarity metric, defined as  $h_{ij}(t):= ||f_i(t)- {f}^{j(i)}_i(t)||$, is the distance between agent $i$'s empirical frequency $f_i(t)$ and the second order belief of agent $i$, i.e., agent $i$'s belief on agent $j$'s belief on $f_i(t)$ denoted with ${f}^{j(i)}_i(t)$. Based on these metrics, agent $i$ decides to communicate its empirical frequency $f_i(t)$ to agent $j$ if the following logical condition is satisfied, 
\begin{align}\label{eq_ack_check}
    \bbone(\eta_1 \le h_{ii}(t) \le \eta_2) \land \bbone(h_{ij}(t) \ge \eta_3)
\end{align}
where $\eta_2> \eta_1\ge 0$ and $\eta_3\ge 0$, $\bbone(\cdot)$ is the indicator function, and $\land$ is the logical AND operator. Condition \eqref{eq_ack_check} determines the set of agents agent $i$ is willing to communicate with at time step $t$, i.e., $\ccalN_i^{out}(t)$. The set of agents that send their empirical frequencies to agent $i$ at time step $t$ is given by $\ccalN_i^{in}(t)= \{j\in\ccalN\setminus \{i\}: i \in \ccalN_j^{out}(t)\}$.

The intuition for the condition in \eqref{eq_ack_check} is as follows. {The novelty metric $h_{ii}(t)$ has to be between a range for agent $i$ to initiate a communication attempt. The novelty $h_{ii}(t)$ is likely to be small when agent $i$ takes the same action for several steps indicating that it may have converged on an action. If $h_{ii}(t)$ is large, it means agent $i$ is undecided, taking a different action from its past set of actions. When $h_{ii}$ is neither too small or too large, agent $i$ attempts to communicate. Agent $i$ only attempts to send its empirical frequency to agent $j$ if it believes agent $j$ does not have an accurate estimate of its empirical frequency, i.e., if $h_{ij}$ is large enough.} 

Given the communication scheme, agent $i$ updates its belief about agent $j$'s empirical frequency at each time step as follows,
\begin{align}\label{eq_info_ex}
f^{i}_{j}(t)=
\begin{cases}
 f_{j}(t), \, \text{if } c_{ji}(t) =1,\\
f^{i}_{j}(t-1), \, \text{otherwise}.
\end{cases}
\end{align}
That is, agent $i$ replaces its estimate on agent $j$'s empirical frequency with the empirical frequency received from agent $j$ upon a successful communication attempt. Otherwise, its estimate remains the same. %\red{Sarper: We did not state lower bound on succesful communication.}

In computing the belief similarity $h_{ij}(t)$, agent $i$ has to form and update beliefs about agent $j$'s belief on its own empirical frequency $f^j_i(t)$. This can be done via an acknowledgement process where each time agent $i$ makes a successful communication attempt to agent $j$, agent $j$ sends back 1-bit acknowledgement signal. We allow the acknowledgement signal to be subject to failures with a Bernoulli variable $b_{ij}(t)\sim \text{Bernoulli}(\beta_{ij}(t))$ with success rate  $0 \le \beta_{ij}(t)\le 1$. We note that the acknowledgement procedure is executed if and only if agent $i$ receives information from agent $j$. Thus, we have $\mathbb{P}(b_{ij}(t)=0| c_{ji}(t)=0) = 1$. Otherwise, we have $\mathbb{P}(b_{ij}(t)= 1 | c_{ji}(t)=1) > \beta_{ij}(t)$. Given the acknowledgement scheme, agent $i$'s second order belief is updated as follows,
\begin{align}\label{eq_first_order_belief}
f^{j(i)}_{i}(t)=
\begin{cases}
 f_{i}(t), \, \text{if } b_{ji}(t) =1,\\
f^{j(i)}_{i}(t-1), \, \text{otherwise}.
\end{cases}
\end{align}
Upon receiving the acknowledgement, agent $i$ knows that its empirical frequency is transmitted to agent $j$, and agent $j$ has updated its belief as per \eqref{eq_info_ex}. In a scenario where $c_{ij}(t)=1$ and $b_{ji}(t)=1$, empirical frequencies and estimates align, i.e., $f^j_i(t)=f^{j(i)}_{i}(t)= f_{i}(t)$.

\begin{remark}
In the information exchange and belief update protocols described above, each agent keeps an estimate of the empirical frequencies of all other agents $\{f^{i}_j(t)\}_{j\in\ccalN}$, an $N\times K$ real-valued matrix, and second order beliefs about other agents' estimates about its empirical frequency $\{f^{j(i)}_i(t)\}_{j\in\ccalN}$, an $N\times K$ real-valued matrix. 
Agent $i$ attempts to send its empirical frequency $f_i(t)$, a real-valued vector of length $K$, to a subset of agents in $\ccalN$ according to the condition in \eqref{eq_ack_check}. In the decentralized FP considered in \cite{swenson2018distributed}, each agent shares their estimates of all the other agents, $\{f^{i}_j(t)\}_{j\in\ccalN}$, an $N\times K$ real-valued matrix, to all of their neighbors at each step.
\end{remark}

\begin{remark}
The condition in \eqref{eq_ack_check} can be less or more selective depending on the constants $\eta_1$, $\eta_2$, and $\eta_3$. The information exchange protocol becomes more selective as $\eta_1$ and $\eta_3$ is increased and $\eta_2$ is decreased close to $\eta_1$. In contrast, if $\eta_1=\eta_3=0$ and $\eta_2$ is large enough, then each agent broadcasts their empirical frequency to all the agents at each step. If $\eta_3=0$, then there is no need for agents to keep second order beliefs as the term with $h_{ij}(t)$ in \eqref{eq_ack_check} is no longer relevant. If $\eta_2$ is large enough and $\eta_3=0$, the information exchange protocol is equivalent to the communication censoring protocol used for distributed optimization algorithms \cite{chen2018ordered,chen2018lag,li2019communication}. 
\end{remark}

\subsection{Limited Information Communication}

Agents share the maximum value and the index of their empirical frequency, i.e., 
\begin{align}
    \upsilon_i(t)&=\max_{k \in \ccalK} f_{ik}^i(t) \label{eq_lim_max}, \\
    \kappa_i(t)&=\argmax_{k \in \ccalK} f_{ik}^i(t) \label{eq_lim_arg},
\end{align}
instead of their empirical frequencies. When an agent $j$ successfully sends the maximum value $\upsilon_j(t)$ and its index $\kappa_j(t)$ \eqref{eq_lim_arg} to agent $i$, agent $i$ needs to reconstruct a well-defined empirical frequency and update its belief $f^i_j(t)$ accordingly. Upon successful communication of $\upsilon_j(t)$ and $\kappa_j(t)$, the reconstructed belief $f^i_j(t)$ has to satisfy 
\begin{align} \label{eq_proper_belief}
    \sum_{k \in \ccalK} f^i_{jk}(t) =1, \;     f^i_{jk}(t) \ge 0, \;
    f^i_{j\kappa_i(t)}(t) \ge \upsilon_i(t), 
\end{align}
where $f^i_{jk}$ denotes the $k$th index. While the first two constraints above define a proper distribution over the space of actions, the third constraint makes sure that the receiving agent uses the information received. There could multiple update rules $\Phi_i(\kappa_j(t), \upsilon_j(t))$ that satisfy the conditions in \eqref{eq_proper_belief}. For instance, one update rule can assume full support on the most frequent action of agent $j$, i.e., $f^i_{j\kappa_j(t)}(t)=1$ and $f^i_{jk}(t)=0$ for $k\in \ccalK\setminus \kappa_j(t)$. Another update rule can assume actions other than the most common are equally likely, i.e., $f^i_{j\kappa_i(t)}(t)=\upsilon_j(t)$ and $f^i_{jk}(t)=(1-\upsilon_j(t))/(|\ccalK|-1)$ for $k\in \ccalK\setminus \kappa_j(t)$.

\begin{remark}
The limited communication protocol described further reduces the information sent per communication attempt to a single real value $\upsilon_i(t)$ and an integer $\kappa_i(t)$.  
\end{remark}
\subsection{Convergence of Communication and Belief Update Protocols}

We describe the specific steps of the DFP  with limited and voluntary communication protocols (DFP-VL) in Algorithm \ref{suboptimal_alg_com}. Step 4 corresponds to the {\it best response} step in Algorithm \ref{suboptimal_alg_DFP}. Steps 5-7 correspond to the information sharing and observation steps in Algorithm \ref{suboptimal_alg_DFP}. Steps 8-9 update the empiricial frequency estimates and second order beliefs.

% \red{The algorithm explained throughout this study is a decentralized fictitious play (DFP) algorithm. In the previous section, we extended DFP with several types of communication protocols such as intermittent and limited information, and voluntary communication. We abbreviate this novel algorithm as DFP-ILV. Algorithm \ref{suboptimal_alg_com} outlines the steps of DFP-ILV. }
%
\begin{algorithm}[H]
   \caption{DFP-VL for Agent $i$}
\label{suboptimal_alg_com}
\begin{algorithmic}[1]
   \STATE {\bfseries Input:} The parameters $\rho,\epsilon, \eta_1, \eta_2$, $\eta_3$.
   \STATE {\bfseries Given:} $f_{-i}^i(0)$, $f_{i}^{-i}(0)$, ${f}_{i}^{j(i)}(0)$  and $a(0)$ for all $i \in \ccalN$.
\FOR{$t=1,2,\cdots $}
  \STATE Agent $i$ takes action $a_i(t)$ using \eqref{eq_response_dc}. 
  \STATE Determine $\ccalN_i^{out}(t)$ by checking \eqref{eq_ack_check} for all $j\in \ccalN \setminus \{i\}$. 
  \STATE Transmit $\upsilon_i(t)$ and $\kappa_i(t)$ to agent $j\in \ccalN_i^{out}$.
  \STATE Send an acknowledgement signal to agent $j\in \ccalN_i^{in}(t)\cap \{j: c_{ji}(t)=1\}$.
  \STATE   Update $\{f^{i}_j(t)\}_{j \in \ccalN \setminus \{i\}}$ using \eqref{eq_info_ex} .
\STATE  For agent $\{j\in \ccalN_i^{out}\cap \{j: b_{ji}(t)=1\}\}$, i.e., if communication and acknowledgement are successful, update $f^{j(i)}_{i}(t) = f_i^i(t)$ accordingly, otherwise $f^{j(i)}_{i}(t) = f^{j(i)}_{i}(t-1)$.
  \ENDFOR 
   \end{algorithmic}
\end{algorithm}

\begin{theorem} \label{lem_prob_com2}
Suppose the communication and acknowledgement success probabilities are lower bounded by a positive value, i.e., $p_{ij}(t)>\nu>0$ and $\beta_{ji}(t)>\nu>0$ for all $t\in \naturals^+$ and $i\in \ccalN, j\in\ccalN$. Let $\{a(t)=(a_1(t),(a_2(t),\cdots,a_N(t))\}_{t\ge1}$ be a sequence of actions generated by the DFP-VL (Algorithm \ref{suboptimal_alg_com}). Then, Condition \ref{cond_prob} is satisfied for any $\xi>0$ given small enough $0\le\eta_1<\xi/2$, large enough $\xi/2< \eta_2$, and small enough $0\le \eta_3\le \xi/2$ such that
if an agent $j\in\ccalN$ repeats the same action for at least $T>\hat T$ times starting from time $t>0$, i.e., $a_{j}(s)=\bbe_k$ for  $s=t,t+1,\cdots,t+T-1$ and $\bbe_k\in\ccalA$, agent $i\in\ccalN$ learns agent $j$'s action with positive  probability $\hat{\epsilon}>0$, i.e., $\mathbb{P} (||a_j(t+T)-f^i_j(t+T)|| \le {\xi}| {\ccalH(t))} \ge \hat{\epsilon}$.
%Condition \ref{cond_prob} is satisfied such that there exists a positive probability that agent $i \in \ccalN$ learns others' static actions repeated for $T> \hat T$ time steps starting from time $t>0$ with positive probability, i.e., $\mathbb{P} (||a_j(t+T)-f^i_j(t+T)|| \le \blue{\xi}| \red{\ccalH(t)}) \ge \hat{\epsilon}$.
\end{theorem}

\begin{myproof}
See Appendix.
\end{myproof}

Theorem \ref{lem_prob_com2} implies that DFP-VL converges to a pure NE of any weakly acyclic game via Theorem \ref{thm_main}. 
% \begin{remark}
% \red{We stated that DFP-IVL algorithm satisfies Condition \ref{cond_prob}. Note that usage of DFP with any combination of three different protocols also assures the convergence to a pure NE. We can still show that they satisfy Condition \ref{cond_prob}, by letting $\eta_1=\sqrt{2}$, $\eta_2=0$, $\eta_3=0$ or, $\phi_j(t+T)=f^j_j(t+T)$.} \green{CE: Rewrite this based on the changes in previous section. Essentially the DFP-LV works with any choice of threshold parameters. }
% \blue{Sarper: Let's talk this again.}
% \end{remark}

% \cite{chen2018ordered,chen2018lag,liu2018coca}. 

%%%%%%%%%%%%%%%%%%%%%%%%%%%%%%%%%%%%%%%%%%%%%%%%%%%%%%%%%%%%%%%%%%%%%%%%%%%
%%%   S E C T I O N %%%%%%%%%%%%%%%%%%%%%%%%%%%%%%%%%%%%%
%%%%%%%%%%%%%%%%%%%%%%%%%%%%%%%%%%%%%%%%%%%%%%%%%%%%%%%%%%%%%%%%%%%%%%%%%%%
%
  %%%%%%%%%%%%%%%%%%%%%%%%%%%%%%%%%%%%%%%%%%%%%%%%%%%%%%%%%%%%%%%%%%%%%%%%%%%
%%%   B E G I N    F I G U R E %%%%%%%%%%%%%%%%%%%%%%%%%%%%%%%%%%%%%
%%%%%%%%%%%%%%%%%%%%%%%%%%%%%%%%%%%%%%%%%%%%%%%%%%%%%%%%%%%%%%%%%%%%%%%%%%%
\begin{figure*}
	\centering
	\begin{tabular}{ccc
	}
	\includegraphics[width=.33\linewidth]{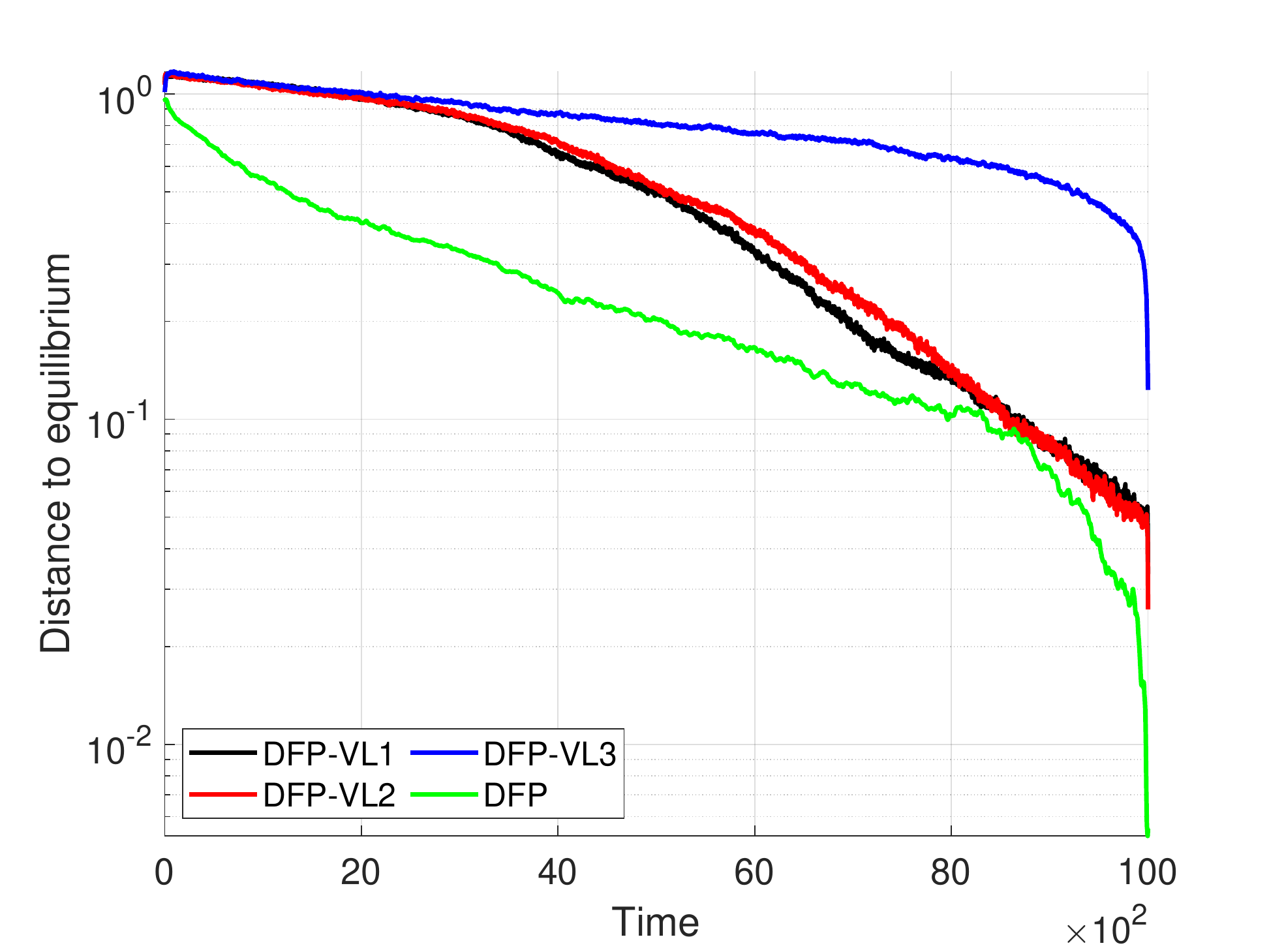}& 
	\includegraphics[width=.33\linewidth]{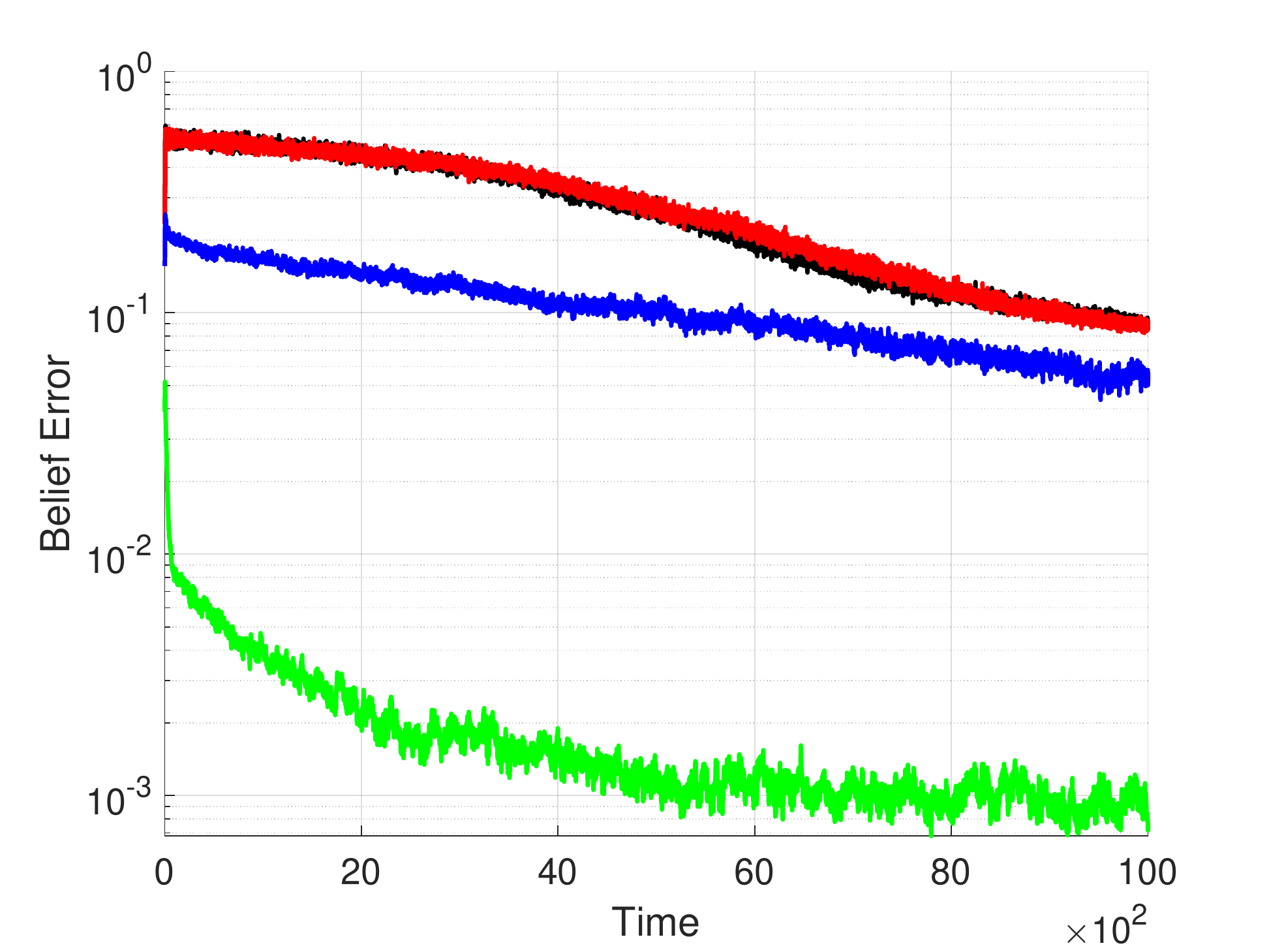} &
    \includegraphics[width=.33\linewidth]{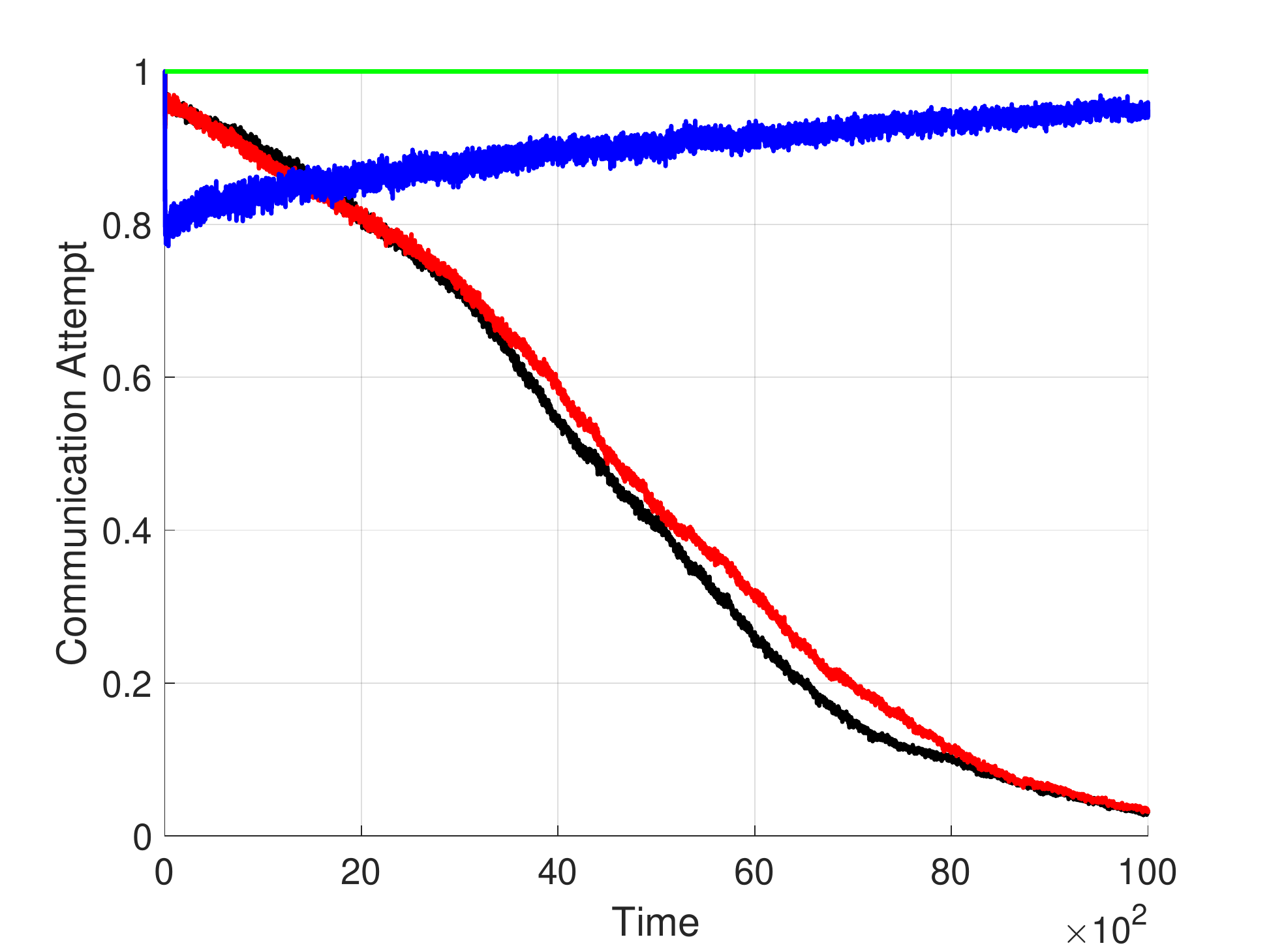} \\ 
	(a) & (b) & (c) \\
    %(c) & (d)
	\end{tabular}
	\caption{Convergence results over 100 replications.  (a) Convergence of empirical frequencies to pure NE $\frac{1}{N}\sum_{i \in \ccalN} || f^i_i(t)-a^*_i||$ on average. We obtain the nearest pure NE by solving a linear assignment problem. (b) Convergence of beliefs $\frac{1}{N(N-1)}\sum_{i \in \ccalN} \sum_{j \in \ccalN \setminus \{i\}}||f^i_i(t)-f^j_i(t) ||$.(c) Average attempt per communication link over time. } \vspace{0pt}
	\label{fig_conv}
\end{figure*}

%(c)	Convergence of Second Order Beliefs $\frac{1}{N(N-1)}\sum_{i \in \ccalN} \sum_{j \in \ccalN \setminus \{i\}}||f^i_i(t)-\hat{f}^j_i(t) ||$. 
%%%%%%%%%%%%%%%%%%%%%%%%%%%%%%%%%%%%%%%%%%%%%%%%%%%%%%%%%%%%%%%%%%%%%%%%%%%
%%%   E N D    F I G U R E %%%%%%%%%%%%%%%%%%%%%%%%%%%%%%%%%%%%%
%%%%%%%%%%%%%%%%%%%%%%%%%%%%%%%%%%%%%%%%%%%%%%%%%%%%%%%%%%%%%%%%%%%%%%%%%%%
%%%%%%%%%%%%%%%%%%%%%%%%%%%%%%%%%%%%%%%%%%%%%%%%%%%%%%%%%%%%%%%%%%%%%%%%%%%

\section{Numerical experiments} \label{sec_numeric}

%To numerically assess the effectiveness of DFP-V, we are going to use Channel Interference Game as the example problem, which will be explained in detail next.

We investigate the performance of different communication protocols in terms of convergence rate and cost of communication in the target assignment game. 

\subsection{Target Assignment Game}

%\red{Target assignment problem is to provide best possible assignment subject to constraints. In multi-agent setting where its recent applications arise in robotic teams, target assignment problem can also be remodeled as game-theoretical problem, which we name it as target assignment game. In this setting,each agent seeks to select an uncovered target with minimum distance. Therefore, the utility function for each agent $i \in \ccalN$ can be given as,  }\green{Remove}

A team of $N$ agents aim to cover $N$ targets with minimum effort. Given the  decentralized decision-making setting, we can represent the problem as a game with the following utility function for agent $i$,

\begin{equation}\label{util_target}
    u_{i}(a_i,a_{-i})= \frac{a_i^T \bbone_{a_{-ik}=0}}{a_i^T d_{i}},
\end{equation}
where $a_i=\bbe_k \in \mathbb{R}^K$ is an unit vector and $\mathbbm{1}_{a_{-ik}=0} \in \{0,1\}^K$ is a binary vector whose $k^{th}$ index is 1 if none of the other agents $j \in \ccalN \setminus \{i\}$ selects target $k$, and otherwise the $k^{th}$ index is equal to 0. 
The distance vector $d_i=[d_{i1},\cdots,d_{ik}, \cdots, d_{iK}] \in \mathbb{R}_+^K $ measures the distance between agent $i$ and each target $k$ in 2D plane, where $d_{ik}=||\theta_i-\theta_k||$. Agent $i$ obtains a positive utility that is inversely proportional to the distance of the agent to the selected target if the  target is not selected by another agent $j$. Otherwise, agent $i$ receives zero utility. Given the utility function  \eqref{util_target}, any joint action that is one-to-one assignment between agents and targets is a pure NE.

In the numerical experiments, we consider {a} target assignment problem with $N=20$ agents and $K=20$ targets. Positions of agents and targets are randomly generated in  a $2$-D plane. 
{Target positions are generated with polar coordinates whose radii are uniformly sampled from $15$ to $20$, and angular coordinates are also uniformly sampled between $0$ and $2\pi$}. Similarly, the positions of agents on the x-y plane are determined by sampling from a normal distribution with mean $0$ and standard deviation 1 independently for each axis. Using the positions, distances between agents and targets are computed.

 %\red{Each coordinate of positions of targets are uniformly sampled around circles whose radius are between $15$ and $20$.} 
%\green{CE: I can't parse this sentence... plural of radius is radii.}

The communication and acknowledgement probability for each link are given as $p_{ij}(t)=0.6$ and $\beta_{ij}(t)=0.9$ for all $t\geq 1$ and pairs of agents. Initial empirical frequencies of agents $f_i(0)$ are set to uniform discrete distribution, i.e., $f_{ik}=1/K$ for $k=\{1, \dots, K\}$. We run each simulation for $T_f=10,000$ steps. 

\subsection{Effects of the Communication Protocol}

We compare three versions of the DFP-VL algorithm to the standard DFP in which agents attempt to communicate with all the agents after each decision. In DFP-VL1 we have all the communication bounds in \eqref{eq_ack_check} relevant. In DFP-VL2, we ignore the upper bound on $h_{ii}(t)$ by making $\eta_2$ large. In DFP-VL3, agents attempt to communicate with all the other agents as long as $h_{ii}(t)$ is bounded by $\eta_2$. See Table~\ref{table:param} for specific parameter values.

\begin{table}[t]
\centering
\begin{tabular}{@{}l l l l l @{}}\toprule
\multicolumn{1}{ c  }{} & \multicolumn{1}{ c  }{}&     \multicolumn{1}{ c  }{Parameters}&            \multicolumn{1}{ c}{}\\
\cmidrule{2-5}
\multicolumn{1}{ c  }{} & \multicolumn{1}{ c }{DFP-VL1} &\multicolumn{1}{ c }{DFP-VL2}   &\multicolumn{1}{ c }{DFP-VL3} &\multicolumn{1}{ c }{DFP} \\
\cmidrule{2-5}
%\multicolumn{1}{ c  }{} & 0.05  & 5  &\\ \midrule
\multicolumn{1}{ c  }{$\eta_1$}& 0.01 & 0.01 &- &- \\
\multicolumn{1}{ c  }{$\eta_2$}& 0.6 & - &0.7  &- \\ 
\multicolumn{1}{ c  }{$\eta_3$}& 0.01 & 0.01 &- &-\\
\multicolumn{1}{ c  }{$\epsilon$}& 0.3 & 0.3 &0.1 &0.9 \\
\multicolumn{1}{ c  }{$\rho$}& 0.6 & 0.6 &0.4 &0.1 \\
\bottomrule
\vspace{0pt}
\end{tabular}
\caption{PARAMETER VALUES OF ALGORITHMS}
\label{table:param}
\end{table}
%\red{(Sarper) Note: here I mention $\eta_3$ as lower bound for $H_{ii}$ (the distance between a pure action and agents' own empirical frequencies.)  }

%We randomly created $100$ instances of targets and positions for each algorithm in order to obtain average behaviour of algorithms. \red{ As the problem size increases, it becomes harder to determine closest pure NE to given final empirical frequencies $f^i_i(T_f)$. This situation requires the development of methodology to find closest pure NE. At the end of each replication, }\green{CE: Not necessary to discuss your hardships. Remove.}\red{we solve linear one-to-one assignment problem to compute average distances to nearest pure NE.}\green{CE: Confused... Why do you say another? }

\begin{comment}
\begin{align}
    \sum_{i \in \ccalN} &\sum_{k \in \ccalK}  |f^i_{ik}(T_f)-a^*_{ik}| \\ 
    &\sum_{i \in \ccalN}a^*_{ik}=1, \: k \in \ccalK, \\
    &\sum_{k \in \ccalN}a^*_{ik}=1 \: i \in \ccalN,\\
    &a^*_{ik} \in \{0,1\},
\end{align}
where $a^*_{ik}$ is the binary variables that provides the action of agent $i$ in closest pure NE $a^*$.
\end{comment}

 DFP achieves fastest convergence rates and ends with closest distances to pure NE on average (Fig.~\ref{fig_conv}(a)). DFP-VL1 and DFP-VL2 achieve comparable convergence rates to DFP. DFP-VL3 is the slowest algorithm but shows convergence at acceptable rate around $0.1$ on average at time $T_f=10,000$. We observe that constant communication in DFP achieves a  faster convergence of beliefs, while voluntary communication have a slower convergence in beliefs as shown in Fig.~\ref{fig_conv}(b). Together, Figs.~\ref{fig_conv}(a-b) signify that communication protocols increase belief error but preserve rate of convergence to an equilibrium.

\begin{figure*}
	\centering
	\begin{tabular}{ccc
	}
	\includegraphics[width=.33\linewidth]{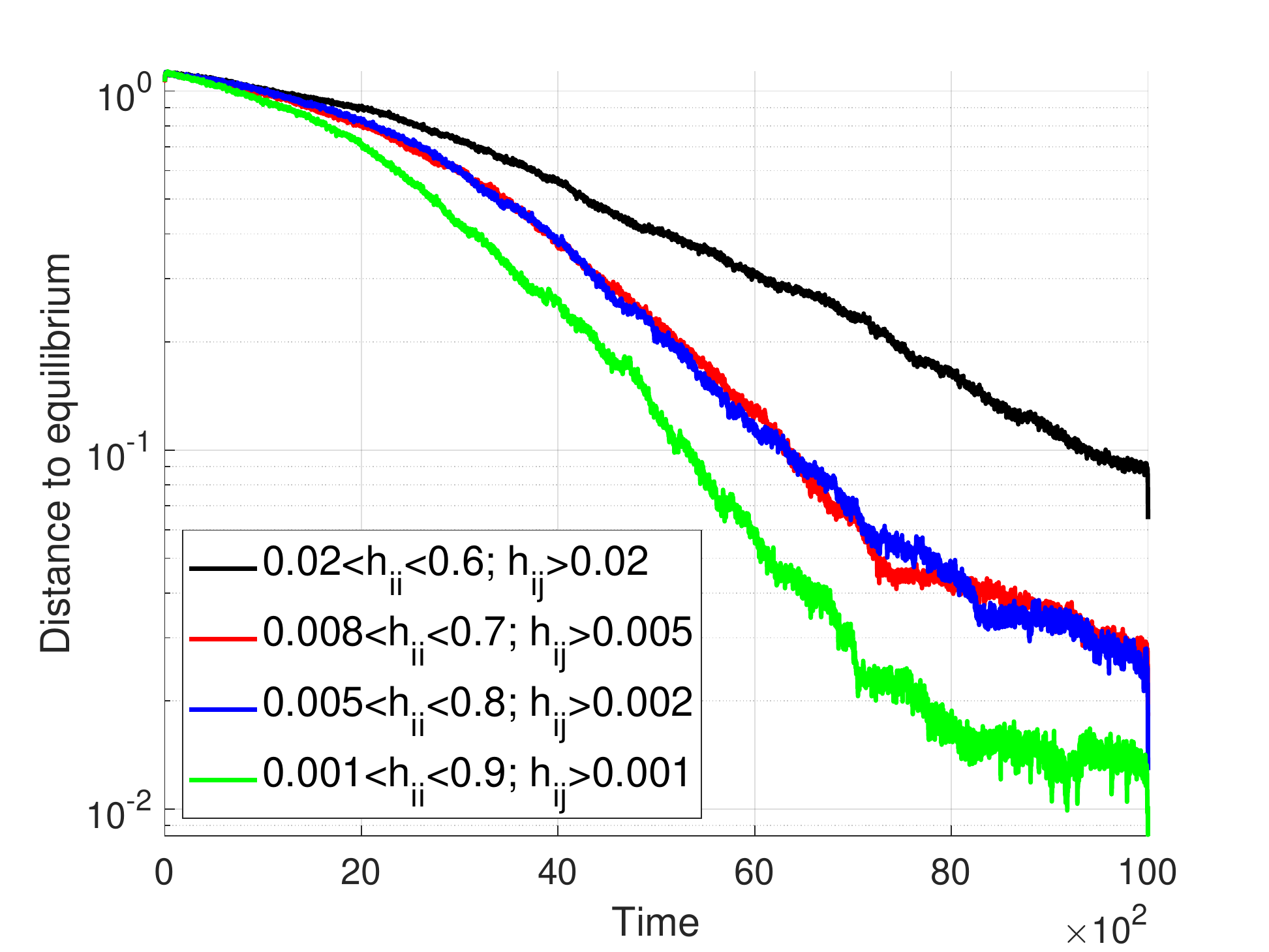}& 
	\includegraphics[width=.33\linewidth]{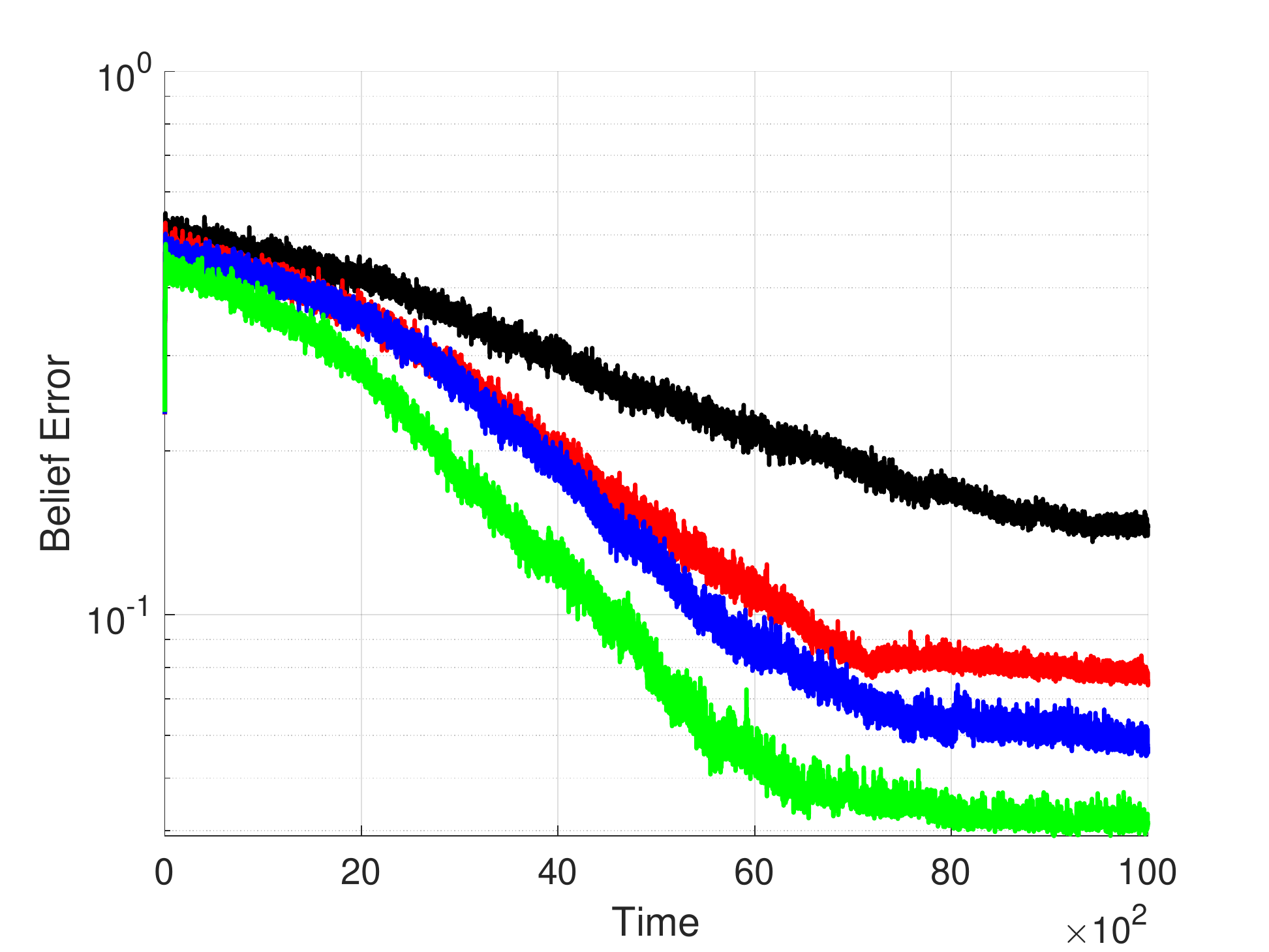}&
	\includegraphics[width=.33\linewidth]{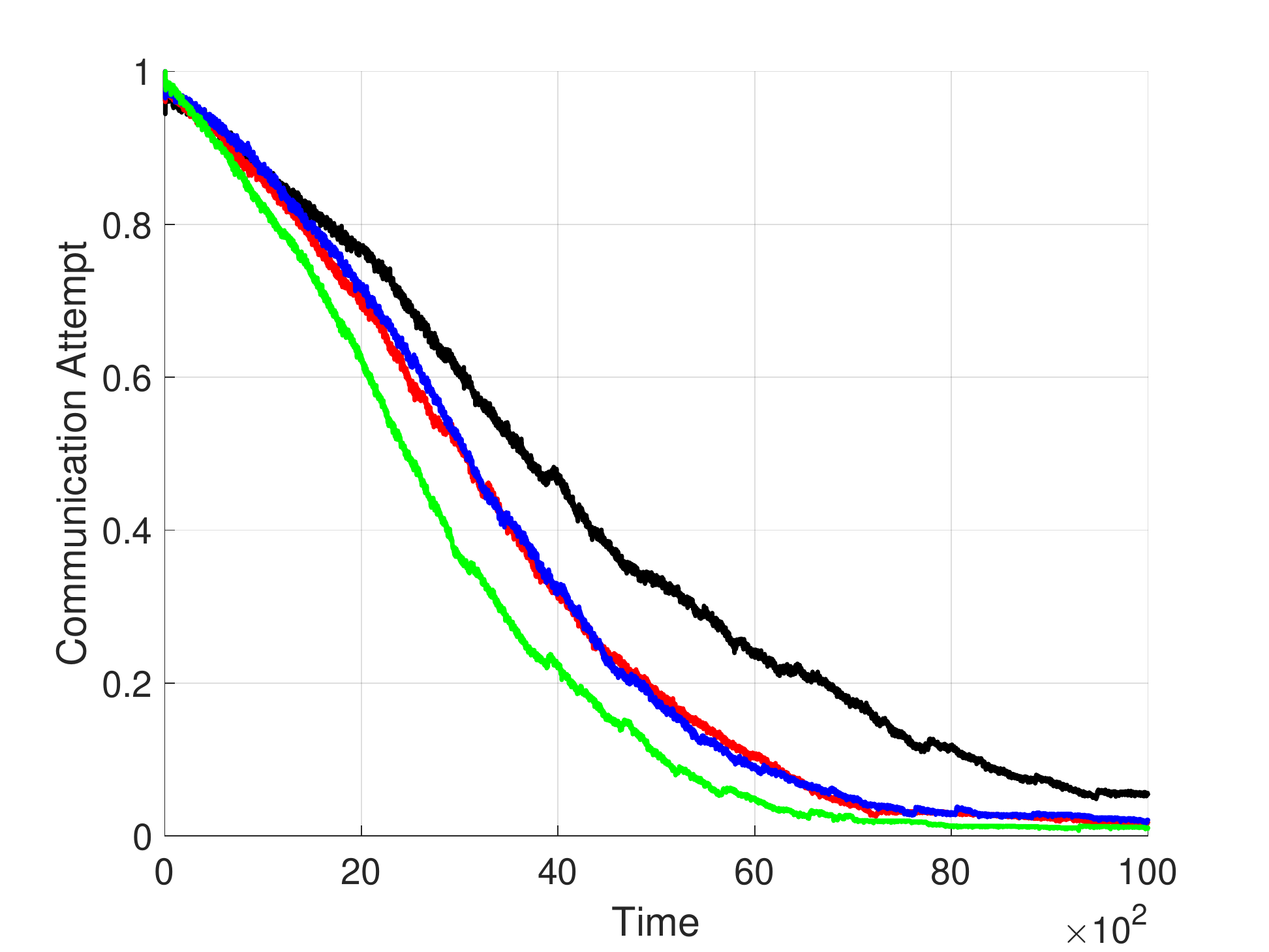} \\
    	(a) & (b) & (c) \\
	\end{tabular}
	\caption{ Convergence results of DFP-VL with different parameter values over $100$ replications. Fading rate $\rho=0.8$ and inertia probability $\epsilon=0.1$. (a) Convergence of empirical frequencies to pure NE $\frac{1}{N}\sum_{i \in \ccalN} || f^i_i(t)-a^*_i||$ on average.\
 (b)	Convergence of Beliefs $\frac{1}{N(N-1)}\sum_{i \in \ccalN} \sum_{j \in \ccalN \setminus \{i\}}||f^i_i(t)-f^j_i(t) ||$.    (c) Average attempt per communication link over time.} \vspace{0pt}
	\label{fig_conv_param}
\end{figure*}
DFP-VL1 and DFP-VL2 utilize $44\%$ and $46\%$ of the communication links, respectively while DFP-VL3 uses $90\%$ of the communication links on average at any point in time (Fig.~\ref{fig_conv}(c)). DFP-VL1 and DFP-VL2 start at full usage of links and then cease the communication attempts almost entirely toward the end of the simulation horizon. Fig.~\ref{fig_conv}(a) and (c) highlight that DFP-VL1 and DFP-VL2 have a faster rate of convergence to NE with a smaller communication effort than DFP-VL3.

\subsection{Parameter Sensitivity}
%\red{ We further expanded numerical results to analyze the effects of parameter tuning. In this part, we focused on DFP-IVL, and used the parameter values of DFP-IVL from the previous set of experiments as a base case and investigated parameter sensitivity. Parameter values can be found at Table-\ref{table:param_sens}.}
%\red{We obtained similar convergence results to pure NE seen at Fig.~\ref{fig_conv}(Top-Left) compared to the previous results DFP-IVL. In addition, there is an improvement over the final convergence results by further reaching around $0.01$. This also implies that DFP-IVL can match the performance of DFP with the proper selection of parameter values.  Fig.~\ref{fig_conv}(Top-Right) and Fig.~\ref{fig_conv}(Bottom-Left) draw parallel conclusions to each other and shows the consistency with convergence rates to pure NE. Smaller error and faster convergence also provides faster converge to pure NE. Again compared to the previous results from DFP-IVL, steeper reduction in communication attempts are visible at Fig.~\ref{fig_conv}(Bottom-Left) and only $27$ percent of possible communication links are attempted over time on average in the best case. For the others, they are still better and $32$ percent for next two ones and $41$ percent of communication links are used in the worst case. 
%\green{Is black line the base case? I could not figure out from the legend. If you do not have the base case in Fig. 2, the reader cannot compare...} 
We assess the performance of DFP-VL under different communication thresholds in Fig. \ref{fig_conv_param}. Here we consider a higher fading and a smaller inertia values compared to Fig. \ref{fig_conv}. Compared with the baseline case (DFP-VL1 shown with black line in Fig. \ref{fig_conv}), we observe that DFP-VL performs better with higher fading and smaller inertia as indicated by the faster convergence to equilibrium by green, blue and red lines in Fig. \ref{fig_conv_param}. The reason is that agents utilize new information faster and update their best-response actions to others accordingly per successful communication when fading is higher and inertia is smaller. In contrast, in standard DFP with constant communication attempts, slow fading rate and increased probability of inertia yields better performance in terms of convergence. This is because in standard DFP, agents are more likely to receive new information at each step. Slow fading rate and increased probability of inertia prevent agents from being oversensitive to new information. We also observe that higher fading and smaller inertia further reduce the communication attempts, e.g., DFP-VL utilizes 27$\%$ of the communication links in Fig. \ref{fig_conv_param}(c) green line. 

Lastly, we observe a counter intuitive phenomenon in Fig. \ref{fig_conv_param}(a-c). As the region of communication is increased, i.e., we decrease $\eta_1$, increase $\eta_2$, and decrease $\eta_3$, not only we get faster convergence and smaller belief errors, but also we see a reduction in communication attempts on average (observe green line with $\eta_1 = 0.001$, $\eta_2 =0.9$ and $\eta_3=0.001$ tread below rest of the lines in Fig. \ref{fig_conv_param}(a-c)). This is counter intuitive because we would expect that a larger region of communication would lead to more communication attempts. However, the frequency of communication attempts is lower in green line except for the first few initial steps. This shows the value of initial communication attempts. Communication in the first few steps allow agents to coordinate early on leading to a smaller belief error, faster convergence to an equilibrium, and reduction in communication attempts later on. Another factor is that smaller $\eta_1$ and $\eta_3$ values yield more precise estimates of agent behavior. This allows agents to be sure early on about the target selections of other agents, and thus eliminate certain targets.

\section{Conclusion}
We considered inertial best-response type algorithms for learning Nash equilibria in weakly acyclic games in random communication networks. We showed that the actions generated from inertial best-response type algorithms converge to a pure Nash equilibrium in weakly acyclic games almost surely under the condition that agents are learn to predict the actions of other agents when those agents repeat the same action. We then proposed voluntary communication protocols for FP in which agents decided who to send the empirical frequencies of their actions based on the novelty of empirical frequency to the receiving agent. We showed that the proposed communication protocols satisfy the prediction under static actions condition, and thus are guaranteed to converge to a pure NE. Compared to standard DFP with constant communication attempts, numerical experiments showed that the voluntary communication protocol significantly reduces communication attempts while achieving a similar convergence rate to a NE.

\appendix

\subsection{Proof of Lemma \ref{lem_est}} \label{ap_pr_2}
%\
\begin{comment}
By Lemma \ref{lem_cond} it holds, $|| f^i_i(t)-\bbe_k||< \xi_1$ and $|| f^j_i(t)-\bbe_k||< \xi_2$ as the result of consecutively taking same action $\bbe_k$ and successful communication attempts by agent $i$ as described. Then, using Assumption \ref{as_max}, there exists a constant $L >0$, $\forall{i} \in \ccalN$ such that, the following holds,
\end{comment}
Since the expectation of the utility function $u_i: \ccalA^N \rightarrow \reals$ is linear and Lipschitz continuous $\forall{i} \in \ccalN$, there exists a constant $L >0$ $\forall{i} \in \ccalN$ such that the following holds, 
\begin{align} \label{eq_lip_up_bound}
    |u_i(a_i,f^i_{-i}(t)-u_i(a_i,a_{-i})|& 
    \le L ||a_{-i}-f^i_{-i}(t) ||,\\
    &=L\sum_{j \in \ccalN \setminus \{i\}}||a_j-f^i_{j}(t)||,\\
    &\le {L (N-1)\xi < \frac{{{\mu}}}{2}}, \label{eq_lipschitz}
\end{align}
for some $\mu>0$. Next, we define the following mutually exclusive subsets of action space $\ccalA$ $\forall{i} \in \ccalN$,
\begin{align}
    \ccalA_1(i)&=\{\bbe_{k_1} \in \ccalA \,| a_i=\bbe_{k_1} \in \argmax u_i(a_i,a_{-i})\}, \\
    \ccalA_2(i)&=\{\bbe_{k_2} \in \ccalA \,| a_i=\bbe_{k_2} \notin \argmax u_i(a_i,a_{-i})\}.
\end{align}
% , suppose that there exist two mutually exclusive subsets of finite action space, $\ccalA$ that holds $\ccalA=\ccalA_1 \cup \ccalA_2$ and $\ccalA_1 \cap \ccalA_2=\emptyset$,  such that they are defined as 
Since they are mutually exclusive subsets, it holds $\ccalA_i=\ccalA_1(i) \cup \ccalA_2(i)$ and $\ccalA_1(i) \cap \ccalA_2(i)=\emptyset$. Then, optimal set over a finite feasible set of utility functions cannot be empty set $\ccalA_1(i)\not =\emptyset$, while it is possible that $\ccalA_2(i)=\emptyset$. Firstly, suppose that $ \ccalA_2(i) \not =\emptyset$. Hence, there exist actions $a'_i \in \ccalA_1(i)$ and $a''_i \in \ccalA_2(i)$ such that,
\begin{equation}\label{eq_opt_dist}
    u_i(a''_i,a_{-i})- \mu < u_i(a'_i,a_{-i}).
\end{equation}
for some $\mu>0$ satisfying \eqref{eq_lip_up_bound}. Note that \eqref{eq_lipschitz} holds for both actions $a_i'\in \ccalA_1(i)$ and $a_i''\in \ccalA_2(i)$,
\begin{align}
    &|u_i(a'_i,f^i_{-i}(t))-u_i(a'_i,a_{-i})| < \frac{\mu}{2}, \label{eq_lipschitz_a} \\
    &|u_i(a''_i,f^i_{-i}(t))-u_i(a''_i,a_{-i})| < \frac{\mu}{2}.\label{eq_lipschitz_b}
\end{align}
Next, we add $u_i(a''_i,f^i_{-i}(t))$ and $u_i(a'_i,f^i_{-i}(t)) $ to the left and right hand sides of \eqref{eq_opt_dist}, respectively. Similarly, we subtract the same corresponding terms from the left and right hand sides of \eqref{eq_opt_dist}. Using the bounds in \eqref{eq_lipschitz_a} and \eqref{eq_lipschitz_b}, we get
\begin{align}
    u_i(a''_i,f^i_{-i}(t))< &u_i(a'_i,f^i_{-i}(t)).
\end{align}
Further, for any two best-response actions, ${a}'_i \in \ccalA_1(i)$ and $\Tilde{a}'_i \in \ccalA_1(i)$, it can be shown that
\begin{align}
    | u_i({a}'_i,f^i_{-i}(t&))- u_i(\Tilde{a}'_i,f^i_{-i}(t))| < {\mu}.
\end{align}
As a result, using its estimates $f^i_{-i}(t)$, agent $i$ only chooses an action from its optimal action set $\ccalA_1(i)$ for the both cases $\ccalA_2(i)= \emptyset$ and $\ccalA_2(i)\not = \emptyset$. Thus, it holds $\forall{i} \in \ccalN$, 
\begin{equation}
    \argmax_{a_i \in \ccalA} u_i(a_i,f^i_{-i}(t))  \subseteq \argmax_{a_i \in \ccalA} u_i(a_i,a_{-i}).
\end{equation}

\subsection{Proof of Theorem \ref{lem_prob_com2}}

We note that the randomness stems from inertia, and communication and acknowledgement failures. The probability of given events in the following part, only depends on these random variables. Thus, showing that the event $\{||a_j(t+T)-f^i_j(t+T||\leq \xi\}$ has a positive probability follows from positive probability of successful communication and acknowledgement, and the positive probability of agent $j$ repeating the same action via inertia. 
%Our proof lies on the idea that using triangle inequality, we firstly can derive an upper bound. Then, we can show that there is positive probability of each event that provides the
% \red{The statement of lemma is to provide the fact that the distance between pure action taken by agent $j$ and the local estimate of agent $i$'s on empirical frequency of agent $j$ has upper bound value $\xi$ with positive probability $\hat{\epsilon}$. Showing the given upper bound holds with positive probability depends on the following events as follows,}
Consider the following events:
\begin{align*}
  E_1(t)=\{&||a_j(t+T)-f^i_j(t+T)|| \le  \\
  &||a_j(t+T)-f_j(t+T)||+||f_j(t+T)-f^{i(j)}_{j}(t)|| \} \\
  E_2(t)=\{&||a_j(t+T)-f_j(t+T)|| \le \xi/2\} \\
  E_3(t)=\{&||f_j(t+T)-f^{i(j)}_{j}(t+T)|| \le \xi/2\}
\end{align*}
\begin{comment}
our aim is to provide lower bound on the given event rather than finding exact value. Therefore, it is enough to show intersection of the following events to have positive probability of $\mathbb{P}(||a_j(t+T)-f^i_j(t+T)|| \le \blue{\xi} | \red{\ccalH(t)}) \ge \hat{\epsilon}$ with given protocols,
\red{ \begin{align} \label{eq_prob_upper}
   \mathbb{P}(&(||a_j(t+T)-f^i_j(t+T)|| \le  \nonumber\\
  &||a_j(t+T)-f^j_j(t+T)||+||f^j_j(t+T)-f^{i(j)}_{j}(t)|| )\cap \nonumber  \\
   &(||a_j(t+T)-f_j(t+T)||+||f^j_j(t+T)-f^{i(j)}_{j}(t+T)|| \nonumber \\
   &\le\xi_1+\xi_2  \le \xi ) \: | \ccalH(t)) \ge \hat{\epsilon}.
\end{align}}
\end{comment}
%\green{CE: I do not follow this conditional probability. First and second events are the same. Make all  $f^j_j$ -> $f_j$}
By triangle equality we have,
\begin{align}
    &||a_j(t+T)-f^i_j(t+T)|| \le  \nonumber\\
    &||a_j(t+T)-f_j(t+T)||+||f_j(t+T)-f^i_j(t+T)||.
\end{align}
Then, via triangle inequality, showing that $E_1(t)$ happens with positive probability reduces to showing the positive probability of the following event,
\begin{align*}
E_4(t)=\{||f_j(t+T)&-f^i_j(t+T)||\le\nonumber \\& ||f_j(t+T)-f^{i(j)}_{j}(t+T)||\}.
\end{align*}
Given the assumptions on $\eta_1$, $\eta_2$ and $\eta_3$, condition \eqref{eq_ack_check} is satisfied, i.e., agent $j$ attempts to communicate with agent $i$, when $E_2$ and $E_3$ is true. 

In the event that agent $j$ successfully communicates with agent $i$ and receives an acknowledgement, we have $f^i_j(t+T)=f^{i(j)}_{j}(t+T)$. Thus, 
\begin{equation}
\mathbb{P}({E_1(t)}|\ccalH(t))=\mathbb{P}({E_4(t)}|\ccalH(t)) \ge  \nu^2,
\end{equation}
where the inequality is via the lower bound on communication and acknowledgement success probabilities. 
% \red{Hence, the probability of the event is at least the product of probability of successful communication and acknowledgement, since in this case, it becomes, $||f_j(t+T)-f^i_j(t+T)||= ||f_j(t+T)-f^{i(j)}_{j}(t)||$,
% \begin{equation}
% \mathbb{P}({E_1(t)}|\ccalH(t))=\mathbb{P}({E_4(t)}|\ccalH(t)) \ge  \red{\epsilon_{com}} \epsilon_{ack}.
% \end{equation}}

Next, the event   $E_2(t)$ is certain given repetition of the same action by agent $j$, and by Lemma \ref{lem_rep} (i) there exists a small enough $\xi_1 \le \xi/2$,

\begin{equation}\label{lem_rep_prob}
    \mathbb{P}(E_2(t) |\ccalH(t)) =1 .
\end{equation}

\begin{comment}
Then, the probability of the event is equal to the event $||\red{f^j_i(t)}-f_j(t)|| \le ||f^{i(j)}_{j}(t)-f_j(t)|| \le \xi_2$, since if the event holds, then it also holds $ ||a_j(t+T)-f^j_j(t+T)||+||f^j_j(t+T)-f^i_j(t+T)|| \le  ||a_j(t+T)-f^j_j(t+T)||+||f^j_j(t+T)-f^{i(j)}_{j}(t)|| \le \xi_1 + \xi_2$. Therefore,
\begin{align}
    \mathbb{P}(&||a_j(t+T)-f^i_j(t+T)|| \le \xi | \ccalH(t) )\\
    = \mathbb{P}(&|||f^j_i(t)-f_j(t)|| \le ||f^{i(j)}_{j}(t)-f_j(t)|| \\ \nonumber
     &\cap ||f^{i(j)}_{j}(t)-f_j(t)|| \le \xi_2  | \ccalH(t) ).
\end{align}
\end{comment}
Now, let $\phi_{j}(t+T)$ be the estimate of empirical frequency of agent $j$ constructed using limited information 
$\upsilon_j(t+T)$ \eqref{eq_lim_max} and $\kappa_j(t+T)$ \eqref{eq_lim_arg} at time $t+T$. By triangle equality, we have
\begin{align} \label{eq_cond2_bound}
    &||f_j(t+T)-f^{i(j)}_{j}(t+T)|| \le ||f_j(t+T)-\phi_{j}(t+T)|| \nonumber\\ &+||\phi_{j}(t+T)-f^i_j(t+T)|| +||f^i_j(t+T)-f^{i(j)}_{j}(t+T)||. 
\end{align}
Now, consider the following events, 
\begin{align*}
  E_5(t)=\{&||f_j(t+T)-\phi_{j}(t+T)|| \le \xi/2\}  \\
  E_6(t)=\{&|||\phi_{j}(t+T)-f^i_j(t+T)|| =0\} \\
  E_7(t)=\{&|||f^i_j(t+T)-f^{i(j)}_{j}(t+T)||=0\}
\end{align*}
Given the repetition of the same actions by agents $j \in \ccalN \setminus \{i\}$ and Lemma \ref{lem_rep} (ii), there exists a small enough $\xi_2 \le \xi/2 $ similar to \eqref{lem_rep_prob}, $\mathbb{P}( E_5(t)| \ccalH(t))=1$.
Further, see the remaining events have also positive probability as the result of successful communication and  acknowledgement, 
\begin{align}
    &\mathbb{P}(E_6(t)| \ccalH(t)) \ge \nu>0, \\
    &\mathbb{P}(E_7(t) | \ccalH(t)) \ge  \nu^2 >0.
\end{align}
\begin{comment}
In the same way, if successful communication and acknowledgement happen, it becomes $f^j_i(t)=f^{i(j)}_{j}(t)$. Hence, the probability of the event is at least the product of successful communication and  acknowledgement,
\begin{align}
    &\mathbb{P}(||f^j_i(t)-f_j(t)|| \le ||f^{i(j)}_{j}(t)-f_j(t)|| | \ccalH(t) ) \nonumber \\
    &\ge  p_{ij}(t)\beta_{ij}(t) \ge \nu^2>0. 
\end{align}
%\blue{Sarper: The lower bound of communication and acknowledgement is $\epsilon$ and $\epsilon_{ack}$. We can reduce to only $\epsilon$.}
\end{comment}
From \eqref{eq_cond2_bound} and the bounds above, we have 
\begin{align}
     \mathbb{P}(E_3(t)| \ccalH(t)) \ge \mathbb{P}(E_5(t),E_6(t),E_7(t)| \ccalH(t))\ge \nu^2 >0.
\end{align}
Thus, there exists a positive real number $\hat{\epsilon}>0$ such that, 
\begin{align}
    &\mathbb{P} (||a_j(t+T)-f^i_j(t+T)|| \le {\xi}| \ccalH(t)) \ge \nonumber\\
    &\mathbb{P}(E_1(t),E_2(t),E_3(t)| \ccalH(t)) \ge \nu^2 =\hat{\epsilon}>0.
\end{align}

\subsection{Technical Result}\label{ap_pr_1}
\begin{lemma}\label{lem_rep}
 Let $\{a(t)=(a_1(t),a_2(t),\cdots,a_N(t))\}_{t\ge1}$ be a sequence of actions  generated by the DFP-VL (Algorithm \ref{suboptimal_alg_com}). Suppose agent $j \in \ccalN \setminus \{i\}$ repeats the same action $a_{j}(s)=\bbe_k$ at least ${T>0}$ times for $s=t,t+1,\cdots,{t+T-1}$. Then, there exist $0 <\xi_1$ and $0 <\xi_2$ such that following  statements hold,
 \begin{itemize}
     \item [\textit{i)} ] $||a_j(t+T)-f_j(t+T)|| \le {\xi_1}$ for all $j\in\ccalN\setminus \{i\}$,
     \item [\textit{ii)} ] $||\phi_j(t+T)-f_j(t+T)|| \le {\xi_2}$ for all $j\in\ccalN\setminus \{i\}$,
 \end{itemize}
 where $\phi_j(t)$ is the reconstructed belief of agent $j$'s empirical frequency using $\upsilon_j(t)$ and $\kappa_j(t)$ defined in \eqref{eq_lim_max} and \eqref{eq_lim_arg}, respectively. 
\end{lemma}
\begin{myproof}
\begin{itemize}
 \item [\textit{i)} ] From \eqref{eq_empirical_frequency}, it holds that if $\bbe_k$ is repeated for any $\tau \in \{0,1,2,\cdots\}$ starting from time $t$ by a player ${j} \in \ccalN \setminus \{i\}$,
\begin{align}
    f_j(t+\tau)=(1-\rho)^{\tau}f_j(t)+(1-(1-\rho)^{\tau})\bbe_k,
\end{align}
    Subtracting $\bbe_k$ from both sides and taking the norm we obtain the following,
\begin{align}\label{eq_rep_1}
    ||f_j(t+\tau)-\bbe_k||&=||(1-\rho)^{\tau}(f_j(t)-\bbe_k)||, \\
    &= O((1-\rho)^{\tau}).
\end{align}
 Therefore, if agent $j \in \ccalN \setminus \{i\}$  repeat the same action $a_{j}(s)=\bbe_k)$  at least {$T>0$} times for {$s=t,t+1,\cdots,t+T-1$}, there exists a positive upper bound ${\xi_1} >0$,
\begin{equation}
    ||a_j(t+T)-f_j(t+T)|| \le {\xi_1}\;\; \forall j\in\ccalN\setminus \{i\}.
\end{equation}
\item [\textit{ii)} ] To provide an upper bound on we can use triangle inequality as below,
\begin{align}
    ||\phi_j(t+T)-f_j(t+T)||& \le  ||\phi_j(t+T)-a_j(t+T)|| \nonumber\\
     &+||a_j(t+T)-f_j(t+T)||.
\end{align}
Then, notice that $||a_j(t+T)-f_j(t+T)||= O((1-\rho)^{T+1})$ implies $|\upsilon_j(t)-1| = O((1-\rho)^{T+1})$. Since $\phi_{j\kappa_i(t)}(t+T)\ge\upsilon_j(t)$ via \eqref{eq_proper_belief}, it also holds $ ||\phi_j(t+T)-a_j(t+T)||= O((1-\rho)^{T+1})$. Thus, given the repetition of the same action,  there exists a positive upper bound ${\xi_2} >0$,
\begin{align}
    &||\phi_j(t+T)-f_j(t+T)||  \nonumber\\
    & \le ||\phi_j(t+T)-a_j(t+T)||+||a_j(t+T)-f_j(t+T)|| \\
    & \le  O((1-\rho)^{T+1}) \le \xi_2.
\end{align}
 \end{itemize}
\end{myproof}

\bibliographystyle{IEEEtran}
\bibliography{bibliography}

% Generated by IEEEtran.bst, version: 1.14 (2015/08/26)
\begin{thebibliography}{10}
\providecommand{\url}[1]{#1}
\csname url@samestyle\endcsname
\providecommand{\newblock}{\relax}
\providecommand{\bibinfo}[2]{#2}
\providecommand{\BIBentrySTDinterwordspacing}{\spaceskip=0pt\relax}
\providecommand{\BIBentryALTinterwordstretchfactor}{4}
\providecommand{\BIBentryALTinterwordspacing}{\spaceskip=\fontdimen2\font plus
\BIBentryALTinterwordstretchfactor\fontdimen3\font minus
  \fontdimen4\font\relax}
\providecommand{\BIBforeignlanguage}[2]{{%
\expandafter\ifx\csname l@#1\endcsname\relax
\typeout{** WARNING: IEEEtran.bst: No hyphenation pattern has been}%
\typeout{** loaded for the language `#1'. Using the pattern for}%
\typeout{** the default language instead.}%
\else
\language=\csname l@#1\endcsname
\fi
#2}}
\providecommand{\BIBdecl}{\relax}
\BIBdecl

\bibitem{eksin2017disease}
C.~Eksin, J.~S. Shamma, and J.~S. Weitz, ``Disease dynamics in a stochastic
  network game: a little empathy goes a long way in averting outbreaks,''
  \emph{Scientific reports}, vol.~7, p. 44122, 2017.

\bibitem{bauch2004vaccination}
C.~T. Bauch and D.~J. Earn, ``Vaccination and the theory of games,''
  \emph{Proceedings of the National Academy of Sciences}, vol. 101, no.~36, pp.
  13\,391--13\,394, 2004.

\bibitem{kar2014distributed}
S.~Kar, G.~Hug, J.~Mohammadi, and J.~M. Moura, ``Distributed state estimation
  and energy management in smart grids: A consensus $+$ innovations approach,''
  \emph{IEEE Journal of selected topics in signal processing}, vol.~8, no.~6,
  pp. 1022--1038, 2014.

\bibitem{zhang2012robust}
Y.~Zhang, N.~Gatsis, and G.~B. Giannakis, ``Robust distributed energy
  management for microgrids with renewables,'' in \emph{2012 IEEE Third
  International Conference on Smart Grid Communications (SmartGridComm)}.\hskip
  1em plus 0.5em minus 0.4em\relax IEEE, 2012, pp. 510--515.

\bibitem{aydin2020communication}
S.~Ayd{\i}n and C.~Eksin, ``Communication censoring in decentralized fictitious
  play for the target assignment problem,'' in \emph{2020 IEEE Conference on
  Control Technology and Applications (CCTA)}.\hskip 1em plus 0.5em minus
  0.4em\relax IEEE, 2020, pp. 334--339.

\bibitem{kantaros2016distributed}
Y.~Kantaros and M.~M. Zavlanos, ``Distributed communication-aware coverage
  control by mobile sensor networks,'' \emph{Automatica}, vol.~63, pp.
  209--220, 2016.

\bibitem{kantaros2019temporal}
Y.~Kantaros, M.~Guo, and M.~M. Zavlanos, ``Temporal logic task planning and
  intermittent connectivity control of mobile robot networks,'' \emph{IEEE
  Transactions on Automatic Control}, 2019.

\bibitem{brown1951iterative}
G.~W. Brown, ``Iterative solution of games by fictitious play,'' \emph{Activity
  analysis of production and allocation}, vol.~13, no.~1, pp. 374--376, 1951.

\bibitem{young2004strategic}
H.~P. Young, \emph{Strategic learning and its limits}.\hskip 1em plus 0.5em
  minus 0.4em\relax OUP Oxford, 2004.

\bibitem{marden2009joint}
J.~R. Marden, G.~Arslan, and J.~S. Shamma, ``Joint strategy fictitious play
  with inertia for potential games,'' \emph{IEEE Transactions on Automatic
  Control}, vol.~54, no.~2, pp. 208--220, 2009.

\bibitem{monderer1996fictitious}
D.~Monderer and L.~S. Shapley, ``Fictitious play property for games with
  identical interests,'' \emph{Journal of economic theory}, vol.~68, no.~1, pp.
  258--265, 1996.

\bibitem{marden2009cooperative}
J.~Marden, G.~Arslan, and J.~Shamma, ``Cooperative control and potential
  games,'' \emph{IEEE Trans. Syst., Man, and Cybern. B, Cybern.}, vol.~39,
  no.~6, pp. 1393--1407, 2009.

\bibitem{robinson1951iterative}
J.~Robinson, ``An iterative method of solving a game,'' \emph{Annals of
  mathematics}, pp. 296--301, 1951.

\bibitem{candogan2013dynamics}
O.~Candogan, A.~Ozdaglar, and P.~A. Parrilo, ``Dynamics in near-potential
  games,'' \emph{Games and Economic Behavior}, vol.~82, pp. 66--90, 2013.

\bibitem{sayin2020fictitious}
M.~O. Sayin, F.~Parise, and A.~Ozdaglar, ``Fictitious play in zero-sum
  stochastic games,'' \emph{arXiv preprint arXiv:2010.04223}, 2020.

\bibitem{Swenson_et_al_2014}
B.~Swenson, S.~Kar, and J.~Xavier, ``Empirical centroid fictitious play: An
  approach for distributed learning in multi-agent games,'' \emph{IEEE Trans.
  Signal Process.}, vol.~63, no.~15, pp. 3888 -- 3901, 2015.

\bibitem{eksin2017distributed}
C.~Eksin and A.~Ribeiro, ``Distributed fictitious play for multiagent systems
  in uncertain environments,'' \emph{IEEE Transactions on Automatic Control},
  vol.~63, no.~4, pp. 1177--1184, 2017.

\bibitem{swenson2018distributed}
B.~Swenson, C.~Eksin, S.~Kar, and A.~Ribeiro, ``Distributed inertial
  best-response dynamics,'' \emph{IEEE Transactions on Automatic Control},
  vol.~63, no.~12, pp. 4294--4300, 2018.

\bibitem{arefizadeh2019distributed}
S.~Arefizadeh and C.~Eksin, ``Distributed fictitious play in potential games
  with time varying communication networks,'' in \emph{2019 53rd Asilomar
  Conference on Signals, Systems, and Computers}.\hskip 1em plus 0.5em minus
  0.4em\relax IEEE, 2019, pp. 1755--1759.

\bibitem{aydin2020decentralized}
S.~Aydin and C.~Eksin, ``Decentralized fictitious play with voluntary
  communication in random communication networks,'' in \emph{2020 59th IEEE
  Conference on Decision and Control (CDC)}.\hskip 1em plus 0.5em minus
  0.4em\relax IEEE, 2020, pp. 337--342.

\bibitem{alpcan2005distributed}
T.~Alpcan and T.~Ba{\c{s}}ar, ``Distributed algorithms for nash equilibria of
  flow control games,'' in \emph{Advances in dynamic games}.\hskip 1em plus
  0.5em minus 0.4em\relax Springer, 2005, pp. 473--498.

\bibitem{koshal2016distributed}
J.~Koshal, A.~Nedi{\'c}, and U.~V. Shanbhag, ``Distributed algorithms for
  aggregative games on graphs,'' \emph{Operations Research}, vol.~64, no.~3,
  pp. 680--704, 2016.

\bibitem{Shamma_Arslan_2005}
J.~Shamma and G.~Arslan, ``Dynamic fictitious play, dynamic gradient play, and
  distributed convergence to nash equilibria,'' \emph{IEEE Trans. Automatic
  Control}, vol.~50, no.~3, pp. 312--327, 2005.

\bibitem{de2019distributed}
C.~De~Persis and S.~Grammatico, ``Distributed averaging integral nash
  equilibrium seeking on networks,'' \emph{Automatica}, vol. 110, p. 108548,
  2019.

\bibitem{scutari2013joint}
G.~Scutari and J.-S. Pang, ``Joint sensing and power allocation in nonconvex
  cognitive radio games: Nash equilibria and distributed algorithms,''
  \emph{IEEE Transactions on Information Theory}, vol.~59, no.~7, pp.
  4626--4661, 2013.

\bibitem{salehisadaghiani2019distributed}
F.~Salehisadaghiani, W.~Shi, and L.~Pavel, ``Distributed nash equilibrium
  seeking under partial-decision information via the alternating direction
  method of multipliers,'' \emph{Automatica}, vol. 103, pp. 27--35, 2019.

\bibitem{ye2021adaptive}
M.~Ye and G.~Hu, ``Adaptive approaches for fully distributed nash equilibrium
  seeking in networked games,'' \emph{Automatica}, vol. 129, p. 109661, 2021.

\bibitem{chen2018ordered}
Y.~Chen, B.~M. Sadler, and R.~S. Blum, ``Ordered transmission for efficient
  wireless autonomy,'' in \emph{2018 52nd Asilomar Conference on Signals,
  Systems, and Computers}.\hskip 1em plus 0.5em minus 0.4em\relax IEEE, 2018,
  pp. 1299--1303.

\bibitem{chen2018lag}
T.~Chen, G.~Giannakis, T.~Sun, and W.~Yin, ``Lag: Lazily aggregated gradient
  for communication-efficient distributed learning,'' in \emph{Advances in
  Neural Information Processing Systems}, 2018, pp. 5050--5060.

\bibitem{li2019communication}
W.~Li, Y.~Liu, Z.~Tian, and Q.~Ling, ``Communication-censored linearized admm
  for decentralized consensus optimization,'' \emph{IEEE Transactions on Signal
  and Information Processing over Networks}, vol.~6, pp. 18--34, 2019.

\bibitem{young1993evolution}
H.~P. Young, ``The evolution of conventions,'' \emph{Econometrica: Journal of
  the Econometric Society}, pp. 57--84, 1993.

\bibitem{milchtaich1996congestion}
I.~Milchtaich, ``Congestion games with player-specific payoff functions,''
  \emph{Games and economic behavior}, vol.~13, no.~1, pp. 111--124, 1996.

\end{thebibliography}

\end{document}